# ON RELATIONS BETWEEN JACOBIANS OF CERTAIN MODULAR CURVES

IMIN CHEN

ABSTRACT. We confirm a conjecture of Merel describing a certain relation between the jacobians of various quotients of $X(p)$ in terms of specific correspondences.

## 1. INTRODUCTION

Let $X = X(p)/\mathbb{Q}$ be the modular curve which classifies elliptic curves with full level $p$-structure. The curve $X$ has an action of $\mathrm{GL}_2(\mathbb{F}_p)$ defined over $\mathbb{Q}$ and for every subgroup $H$ of $G$, there is a quotient $\pi_H : X \to X_H$ which is also defined over $\mathbb{Q}$.

The subject of this paper concerns a relation between the jacobians of $X_H$ for certain subgroups $H$ of $G$. This relation was verified in [4] [5] using the trace formula, and later by Edixhoven [10] using the representation theory of $G$. It has been noted in way or another by several people including Gross [15], Ligozat [21], Elkies [7].

To describe this relation, suppose now that $p$ is an odd prime and denote by $J_H$ the jacobian of the quotient curve $X_H$. The group $G$ also acts on $\mathbb{P}^1(\mathbb{F}_{p^2}) = \mathbb{P}^1(\mathbb{F}_p[\sqrt{\lambda}])$, where $\left(\frac{\lambda}{p}\right) = -1$, from which we define subgroups $B$, $N$, $N'$ as the stabilisers in $G$ of $\infty$, $\{\infty, 0\}$, $\{\sqrt{\lambda}, -\sqrt{\lambda}\}$, respectively. The relation of jacobians which concerns us is then

**Theorem 1** (Chen,Edixhoven).

$$J_{N'} \times J_B \text{ is isogenous over } \mathbb{Q} \text{ to } J_N \times J_G$$

where we have included $J_G$ (which is trivial in this case) to indicate the form of the relation in contexts other than jacobians.

The results in [4] [5] and [10] only show the existence of the relation of jacobians in Theorem 1 and leave open the question of describing this isogeny explicitly. Subsequently, Merel conjectured in [8] [24] an explicit description of this relation in terms of certain natural correspondences.

This paper will confirm Merel's conjectural description (Theorem 2). To explain this description, we introduce some terminology. A sequence

$$\ldots \longrightarrow A_{i-1} \xrightarrow{\phi_{i-1}} A_i \xrightarrow{\phi_i} A_{i+1} \longrightarrow \ldots$$

will be called almost-exact if the kernel of each morphism $\phi_i$ contains the image of its predecessor $\phi_{i-1}$ with finite quotient $\Phi(A_i)$. Using this terminology, Theorem 1 can be rephrased as

*Date*: 1 August 1998.
1991 *Mathematics Subject Classification.* Primary 11G18, Secondary 14H35.
This research was supported by an NSERC postdoctoral fellowship and a grant from CICMA.

**Theorem 1** (reformulation). *There exist homomorphisms of abelian varieties $(\phi_{GB}, \phi_{BN}, \phi_{NN'})$ defined over $\mathbb{Q}$ such that*

$$0 \longrightarrow J_G \xrightarrow{\phi_{GB}} J_B \xrightarrow{\phi_{BN}} J_N \xrightarrow{\phi_{NN'}} J_{N'} \longrightarrow 0 \tag{1}$$

*is almost-exact.*

For $H \subset K$ subgroups of $G$, the universal property of quotients gives a unique map $\pi_{H/K} : X_K \to X_H$. Given a quotient $\pi_H : X \to X_H$, we denote by $\pi_H^* : J_H \to J$ and $\pi_{H*} : J \to J_H$ the homomorphisms of abelian varieties which are induced by Picard and Albanese functoriality. Since $G$ acts on $X$, there is an action of $G$ on $J$ by Albanese functoriality so that an element of $\mathbb{Z}[G]$ yields an endomorphism of $J$ defined over $\mathbb{Q}$.

Let $C$ and $C'$ are the stabilisers in $G$ of $(\infty, 0)$ and $(\sqrt{\lambda}, -\sqrt{\lambda})$, respectively. Then $N = C \cup \omega C$ and $N' = C' \cup \omega' C'$ are the normalisers in $G$ of $C$ and $C'$, respectively, where

$$\omega = \begin{pmatrix} 0 & 1 \\ 1 & 0 \end{pmatrix}$$

$$\omega' = \begin{pmatrix} 1 & 0 \\ 0 & -1 \end{pmatrix}.$$

With the above comments in mind, the conjectural description of the relation of jacobians in Theorem 1 alluded to earlier [8] [24] which will be proved in this paper is the following.

**Theorem 2** (Merel's conjecture). *A choice of $(\phi_{GB}, \phi_{BN}, \phi_{NN'})$ in Theorem 1 is given by*

$$(\phi_{GB}, \phi_{BN}, \phi_{NN'}) \tag{2}$$

$$= (\pi_{G*} \circ \pi_B^*, \pi_{N*} \circ |G|\,(1 - \mathrm{pr}_G)(1 + \omega) \circ \pi_B^*, \left[\frac{p+1}{2}\right] \circ \pi_{N'*} \circ \pi_N^*).$$

*where for a subgroup $H$ of $G$, we let $\mathrm{pr}_H = \frac{1}{|H|} \sum_{h \in H} h \in \mathbb{Q}[G]$ be the projector to the $H$-invariants.*

We note that these choices of $(\phi_{GB}, \phi_{BN}, \phi_{NN'})$ are not optimal. For instance, we may replace $\phi_{NN'}$ by the simpler homomorphism $\pi_{N'*} \circ \pi_N^*$ and still retain almost-exactness. We have chosen these particular homomorphisms, which are multiples of the ones conjectured by Merel, because they are the choices which are induced by natural operators coming from representation theory (see Lemma 6.5).

The method of proof involves reducing the almost-exactness of the above sequence of jacobians to the exactness of an associated sequence of $\mathbb{C}[G]$-modules, where the $\mathbb{C}[G]$-module homomorphisms in the sequence are given by certain double coset operators. An easy application of character theory shows that this sequence of $\mathbb{C}[G]$-modules is exact if and only if one can show certain character sums involving the number of points on a non-trivial family of elliptic curves modulo $p$ are non-zero. Although these character sums seem rather intractable, we exhibit some relations between double coset operators which allow one to simplify them to Legendre character sums [11] [27] and Soto-Andrade sums [18] [30]. Having done so, we show that they are non-zero by using the fact that $p$ does not ramify in the cyclotomic fields generated by these character sums.

The appearance of Legendre character sums and Soto-Andrade sums is suggestive. It is known that natural descriptions of the representations of $\mathrm{SL}_2(\mathbb{F}_p)$ induced from Borel subgroups are described in terms of finite field analogues of hypergeometric functions [13] [14], following indications in [16] [20]. These finite field hypergeometric functions are in turn related to Legendre character sums [27] in a manner similar to the relation between hypergeometric functions and Legendre polynomials in the classical case. Similarly, Soto-Andrade sums appear in the analysis of the representations which are induced from non-split Cartan subgroups [30].

The relation of jacobians in Theorem 1 has connections to a question of Serre [28] which asks whether there are only finitely-many primes $p$ for which there exists a non-CM elliptic curve over $\mathbb{Q}$ whose the mod $p$ representation is non-surjective. Indeed, if the mod $p$ representation of an elliptic curve over $\mathbb{Q}$ is non-surjective, then it has image lying in a maximal proper subgroup of $G$, i.e. a conjugate of $B$, $N$, $N'$, or an exceptional group whose projective image is $S_4$ (the groups $A_4$ and $A_5$ do not occur due to properties of the Weil pairing). However, such elliptic curves are essentially classified by the quotient curves $X_B$, $X_N$, $X_{N'}$, $X_{S_4}$ [9] so that such an elliptic curve gives rise to a $\mathbb{Q}$-rational point on one of these curves. The question can therefore be rephrased as whether there are only finitely-many primes $p$ such that the quotient curves $X_B$, $X_N$, $X_{N'}$, $X_{S_4}$ have a non-cuspidal, non-CM $\mathbb{Q}$-rational point.

Mazur has answered the question in the affirmative for the modular curve $X_B \cong_{\mathbb{Q}} X_0(p)$ [23]. Indeed it is shown that for $p > 163$, the only $\mathbb{Q}$-rational points on $X_0(p)$ are cuspidal. In [22, p. 36], it is noted that Serre has shown that $X_{S_4}$ has no $\mathbb{Q}$-rational points for $p > 13$. Momose has obtained some partial results for the modular curve $X_N$ [25], in particular when $J_0(p)^-$ has rank 0.

For the quotient curve $X_{N'}$, the relation of jacobians in Theorem 1 suggests that obtaining information about the $\mathbb{Q}$-rational points of the curve $X_{N'}$ via its jacobian is difficult [5] since by standard conjectures about $L$-functions of abelian varieties over $\mathbb{Q}$ [31], this abelian variety does not have any non-trivial quotients with finite Mordell-Weil group over $\mathbb{Q}$.

Serre's question enters into the analysis of diophantine problems similar to Fermat's Last Theorem [24], [26], [8]. However, in these contexts, the elliptic curves in question have extra level structure. This fact was used in [8] to give a solution of Dénes' conjecture by considering a variant of the relation of jacobians in Theorem 1 [7].

We remark that Edixhoven's proof [10] indicates a general phenomenon that relations between jacobians of quotients of a curve arise from relations between idempotents in the group ring associated to the automorphism group of the curve. Indeed, as a general fact, this was already recorded in [17] following work of Accola [1] [2]. We note that an equivalent formulation in terms of character identities is also given in [17].

In the case of Theorem 1, the relevant idempotent relation which was proved by Edixhoven [10] using the character table of $G$ is

**Theorem 3** (Edixhoven).

$$(\mathrm{pr}_{N'} + \mathrm{pr}_B)(1 - \mathrm{pr}_G) \text{ is } \mathbb{Q}[G]^\times\text{-conjugate to } \mathrm{pr}_N(1 - \mathrm{pr}_G).$$

This is equivalent using Frobenius reciprocity to the character identity

**Theorem 4.**

$$1_{N'} + 1_B = 1_N + 1_G,$$

where $1_H$ denotes the character of the permutation representation $\mathbb{Q}[G/H]$ associated to $H$ (i.e. the trivial representation over $\mathbb{Q}$ of $H$ induced to $G$). We note that Arsenault has subsequently given a direct calculation of the above character identity in his thesis [3].

Although it is not too difficult to show the existence of the above idempotent relation via character theory, Merel's conjecture is in some sense a step towards understanding why this relation exists in an explicit way.

## 2. Acknowledgements

I would like to thank the Department of Mathematics at UC Berkeley, the Department of Mathematics and Statistics at McGill University, and l'Institut de Recherche Mathématique de Rennes for their hospitality during my stay there in 1996-97, 1997-98, and 1998, respectively.

Special thanks are due to L. Merel for rekindling my interest on this topic, and whose interest and many stimulating discussions motivated me. I would also like to thank M. Baker, H. Darmon, S. Edixhoven, and N. Yui for many useful discussions and suggestions.

## 3. Reduction to representation theory

In this section, we show how the almost-exactness of a sequence of jacobians such as the one in Theorem 1 and 2 can be deduced from the almost-exactness of an associated sequence of $\mathbb{Z}[G]$-modules.

The action of $G$ on the abelian group $J(\mathbb{C})$ by Albanese functoriality allows one to regard $J(\mathbb{C})$ as a $\mathbb{Z}[G]$-module. For any subgroup $K$ of $G$, we then have the identifications $J(\mathbb{C})^K = \text{Hom}_{\mathbb{Z}[K]}(1, J) = \text{Hom}_{\mathbb{Z}[G]}(\mathbb{Z}[G/K], J(\mathbb{C}))$ by Frobenius reciprocity.

Let $C$ be a curve defined over $\mathbb{Q}$. Its jacobian is also defined over $\mathbb{Q}$. Denote by $\text{Pic}^0(C)(\overline{\mathbb{Q}})$ the $\mathbb{Z}$-module of divisors modulo linear equivalence which are defined over $\overline{\mathbb{Q}}$. By Abel-Jacobi, $\text{Pic}^0(C)(\overline{\mathbb{Q}})$ can be identified with $J(C)(\overline{\mathbb{Q}})$ in a way which is compatible with the action of $G_{\mathbb{Q}}$ [19] [32]. The multiplication by $[n]$ homomorphism on $J$ is described on divisors by sending a divisor class $D$ to $nD$. In particular, since $[n]$ is an isogeny on $J$, it follows for every divisor class $D$, there exists a divisor class $D'$ such that $nD' \sim D$.

Let $X, Y$ be curves defined over $\mathbb{Q}$ and suppose $\pi : X \to Y$ is a non-constant map, also defined over $\mathbb{Q}$. The homomorphisms of jacobians $\pi^* : J(Y) \to J(X)$ and $\pi_* : J(X) \to J(Y)$ induced by Picard and Albanese functoriality, respectively, are defined over $\mathbb{Q}$. Moreover, on the level of divisors, they can be described as follows. The homomorphism $\pi^*$ sends a divisor $D = \sum_{y \in Y} \alpha_y y$ on $Y$ to its pullback divisor $\pi^*(D) = \sum_{y \in Y} \sum_{x \in \pi^{-1}(y)} \alpha_y e_{\pi,x} x$ on $X$ where $e_{\pi,x}$ denotes the ramification index of $x$ with respect to the map $\pi$. The homomorphism $\pi_*$ sends a divisor $D = \sum_{x \in X} \alpha_x x$ on $X$ to its push-down divisor $\pi_*(D) = \sum_{x \in X} \alpha_x \pi(x)$. The above implies for instance that $\pi_* \circ \pi^* = [\deg(\pi)]$. These descriptions follow from the theory of correspondences [29] [32] [19].

We begin with two lemmas which concern the subgroup variety of J on which $H$ acts trivially and the image of the jacobian $J_H$ under the morphism $\pi_H^*$.

**Lemma 3.1.** *Let $J^H = \cap_{h \in H} \ker(h - 1)$. Let $\Sigma_H : J \to J$ be the morphism $\Sigma_H : x \mapsto \sum_{h \in H} hx$. Then the connected component of the identity, $J^{H^o}$, is an abelian subvariety of $J$ which is equal to $\operatorname{im} \sigma_H$, and the group of connected components of $J^H$ is killed by $|H|$.*

*Proof.* The image $\operatorname{im} \Sigma_H$ is an abelian subvariety of $J$. On $\mathbb{C}$-points, we have that $J^H(\mathbb{C}) = J(\mathbb{C})^H$, and $\operatorname{im} \Sigma_H \subset J^{H^o}(\mathbb{C})$ (since it lies in $J^H(\mathbb{C})$, is connected, and contains the identity). Furthermore, the morphism $\Sigma_H$ when restricted to $J^H(\mathbb{C})$ is multiplication by $|H|$ so the image $\operatorname{im} \Sigma_H$ contains $J^{H^o}(\mathbb{C})$. This also shows that $|H|$ kills the group of connected components of $J^H$. $\square$

**Lemma 3.2.** *Let $\pi_H : X \to X_H$ be the quotient map. Then $\operatorname{im} \pi_H^*$ is contained in $J^H$ and the morphisms*

$$\pi_H^* : J_H \to J^{H^o} \subset J$$
$$\pi_{H*} : J^{H^o} \subset J \to J_H$$

*are isogenies with kernel killed by $|H|$.*

*Proof.* For $h \in H$, we see that $\pi_H \circ h = \pi_H$. Hence, by Picard functoriality, $h^* \circ \pi_H^* = \pi_H^*$. Since $h_* \circ h^* = [\deg(h)] = 1$, it follows that $\pi_H^* = h_* \circ \pi_H^*$ so that $H$ acts trivially on $\operatorname{im} \pi_H^*$ and hence $\operatorname{im} \pi_H^* \subset J^{H^o}$.

We therefore have the following sequence of morphisms.

$$J_H \xrightarrow{\pi_H^*} J^{H^o} \xrightarrow{\pi_{H*}} J_H$$
$$J^{H^o} \xrightarrow{\pi_{H*}} J_H \xrightarrow{\pi_H^*} J^{H^o}$$

Since $\pi_{H*} \circ \pi_H^* = [\deg(\pi_H)] = [|H|]$, we see that $\pi_H^*$ is surjective and $\pi_{H*}$ has finite kernel. The morphism $\pi_H^* \circ \pi_{H*}$ when restricted to $J^{H^o}$ is also multiplication by $|H|$ so that $\pi_{H*}$ is surjective and $\pi_H^*$ has finite kernel. Indeed, given a point $P \in J^{H^o}(\mathbb{C}) \subset J(\mathbb{C})^H$, $P$ corresponds to a divisor class $D$ which is invariant under the action of $H$. Thus, $\pi^* \circ \pi_*(D) \sim |H| D$ as desired. $\square$

The starting point for reducing Theorem 1 and 2 to a question in representation theory stems from the following lemma.

**Lemma 3.3.** *Let $\sigma : \mathbb{Z}[G/H'] \to \mathbb{Z}[G/H]$ be a $\mathbb{Z}[G]$-module homomorphism. Then $\sigma$ induces a homomorphism of abelian varieties $\phi : J_H \to J_{H'}$ which is defined over $\mathbb{Q}$.*

*Proof.* By Lemma 3.2, $H$ acts trivially on $\operatorname{im} \pi_H^*$. Thus, $\sigma(H')$ gives a well-defined morphism $\operatorname{im} \pi_H^* \to J$. The desired morphism is then given by $\phi = \pi_{H'*} \circ \sigma(H') \circ \pi_H^*$. $\square$

Indeed, it will be shown in Proposition 3.7 that an almost-exact sequence of $\mathbb{Z}[G]$-modules

$$\ldots \longleftarrow \mathbb{Z}[G/H_{i-1}] \xleftarrow{\sigma_{i-1}} \mathbb{Z}[G/H_i] \xleftarrow{\sigma_i} \mathbb{Z}[G/H_{i+1}] \longleftarrow \ldots$$

induces an almost-exact sequence of jacobians

$$\ldots \longrightarrow J_{H_{i-1}} \xrightarrow{\phi_{i-1}} J_{H_i} \xrightarrow{\phi_i} J_{H_{i+1}} \longrightarrow \ldots$$

Note that this immediately gives Theorem 1 since the character relation in Proposition 4 implies there exists an exact sequence of $\mathbb{Q}[G]$-modules

$$0 \longleftarrow \mathbb{Q}[G/G] \longleftarrow \mathbb{Q}[G/B] \longleftarrow \mathbb{Q}[G/N] \longleftarrow \mathbb{Q}[G/N'] \longleftarrow 0$$

and hence an almost-exact sequence of $\mathbb{Z}[G]$-modules

(3)
$$0 \longleftarrow \mathbb{Z}[G/G] \longleftarrow \mathbb{Z}[G/B] \longleftarrow \mathbb{Z}[G/N] \longleftarrow \mathbb{Z}[G/N'] \longleftarrow 0.$$

To prove Proposition 3.7, we first establish some lemmas.

**Lemma 3.4.** *Suppose we have an almost-exact sequence of $\mathbb{Z}[G]$-modules*

(4)
$$\ldots \longleftarrow M_{i-1} \xleftarrow{\sigma_{i-1}} M_i \xleftarrow{\sigma_i} M_{i+1} \longleftarrow \ldots$$

*For any $\mathbb{Z}[G]$-module $M$, the sequence obtained by taking $\mathrm{Hom}_{\mathbb{Z}[G]}(\cdot, M)$ is almost-exact. Furthermore, $\#\Phi(\mathrm{Hom}_{\mathbb{Z}[G]}(M_i, M)) \leq \#\Phi(M_i) \cdot \#\eta(M_{i-1}, \mathrm{im}\,\sigma_{i-1})$, where $\eta(V, W)$ for $W \subset V$ $\mathbb{Z}[G]$-modules is a quantity to be explained below.*

*Proof.* This essentially follows from the semi-simplicity of $\mathbb{Q}[G]$-modules. Consider the sequence induced by $\mathrm{Hom}_{\mathbb{Z}[G]}(\cdot, M)$

(5) $\quad \mathrm{Hom}_{\mathbb{Z}[G]}(M_{i-1}, M) \xrightarrow{\phi_{i-1}} \mathrm{Hom}_{\mathbb{Z}[G]}(M_i, M) \xrightarrow{\phi_i} \mathrm{Hom}_{\mathbb{Z}[G]}(M_{i+1}, M)$

Since $\sigma_{i-1} \circ \sigma_i = 0$, it follows that $\phi_i \circ \phi_{i-1} = 0$. Hence, $\ker \phi_i \supset \mathrm{im}\,\phi_{i-1}$.

Suppose $W \subset V$ are $\mathbb{Z}[G]$-modules. By semi-simplicity of $\mathbb{Q}[G]$-modules, $V \otimes \mathbb{Q} = W \otimes \mathbb{Q} \oplus W'$ for some complementary $\mathbb{Q}[G]$-module $W'$. Put $W^\perp = W' \cap V$. Let $\eta(V, W)$ be the quotient $\frac{V}{W \oplus W^\perp}$ and $\mathrm{pr}_W$ the projection to $W$ from $W \oplus W^\perp$.

Put $n = \#\Phi(M_i)$, $m = \#\eta(M_{i-1}, \mathrm{im}\,\sigma_{i-1})$ and let $f \in \ker \phi_i$. We claim that $mn \cdot f \in \mathrm{im}\,\phi_{i-1}$. Since $f \in \ker \phi_i$, $n \cdot f \in \ker \phi_i$ so that $\mathrm{im}\,\sigma_i \subset \ker n \cdot f$. Since $n$ kills $\Phi(M_i)$, it follows that in fact $\ker \sigma_i \subset \ker n \cdot f$. Therefore, there exists a $\mathbb{Z}[G]$-module homomorphism $g : \mathrm{im}\,\sigma_{i-1} \to M$ such that $n \cdot f = g \circ \sigma_{i-1}$. One can extend $g$ to a $\mathbb{Z}[G]$-module homomorphism $\tilde{g}$ on all of $M_{i-1}$ by defining $\tilde{g}(x) = g \circ \mathrm{pr}_{\mathrm{im}\,\sigma_{i-1}}(mx)$. Moreover, we see that $\tilde{g} \circ \sigma_{i-1} = mn \cdot f$ so that $\phi_{i-1}(\tilde{g}) = mn \cdot f$. □

**Lemma 3.5.** *Let $\sigma : \mathbb{Z}[G/H'] \to \mathbb{Z}[G/H]$ be a $\mathbb{Z}[G]$-module homomorphism. Let $\phi' : J(\mathbb{C})^H \to J(\mathbb{C})^{H'}$ be the $\mathbb{Z}[G]$-module homomorphism induced by $\sigma$ via the identification $J(\mathbb{C})^K = \mathrm{Hom}_{\mathbb{Z}[G]}(\mathbb{Z}[G/K], J(\mathbb{C}))$. Then we have a commutative diagram*

(6)
$$\begin{array}{ccc} J_H(\mathbb{C}) & \xrightarrow{\phi} & J_{H'}(\mathbb{C}) \\ \pi_H^* \downarrow & & \uparrow \pi_{H*} \\ J(\mathbb{C})^H & \xrightarrow{\phi'} & J(\mathbb{C})^{H'} \end{array}$$

*where $\phi$ is the morphism induced by $\sigma$ in the previous lemma.*

*Proof.* It suffices to show as elements of $\mathrm{Hom}_{\mathbb{Z}[G]}(J(\mathbb{C})^H, J(\mathbb{C})^{H'})$, we have the equality $\phi' = \sigma(H')$. The correspondence between $J(\mathbb{C})^H$ and $\mathrm{Hom}_{\mathbb{Z}[G]}(\mathbb{Z}[G/H], J(\mathbb{C}))$ is such that $x \mapsto (f : g \mapsto gx)$. The $\mathbb{Z}[G]$-module homomorphism which $\sigma$ induces sends $f$ to $f \circ \sigma$. Finally, the correspondence between $\mathrm{Hom}_{\mathbb{Z}[G]}(\mathbb{Z}[G/H'], J(\mathbb{C}))$ and $J(\mathbb{C})^{H'}$ is such that $f \circ \sigma$ is sent to $f \circ \sigma(H') = \sigma(H')x$. □

**Lemma 3.6.** *A sequence of abelian varieties over $\mathbb{C}$*

(7) $$\ldots \longrightarrow J_{i-1} \longrightarrow J_i \longrightarrow J_{i+1} \longrightarrow \ldots$$

*is almost-exact if and only if the associated sequence of $\mathbb{Z}$-modules*

(8) $$\ldots \longrightarrow J_{i-1}(\mathbb{C}) \longrightarrow J_i(\mathbb{C}) \longrightarrow J_{i+1}(\mathbb{C}) \longrightarrow \ldots$$

*is almost-exact. Furthermore, $\#\Phi(J_i)(\mathbb{C}) = \#\Phi(J_i(\mathbb{C}))$.*

*Proof.* Morphisms between varieties over $\mathbb{C}$ are determined uniquely on $\mathbb{C}$-points. □

**Proposition 3.7.** *Suppose the sequence of $\mathbb{Z}[G]$-modules*

$$\ldots \longleftarrow \mathbb{Z}[G/H_{i-1}] \xleftarrow{\sigma_{i-1}} \mathbb{Z}[G/H_i] \xleftarrow{\sigma_i} \mathbb{Z}[G/H_{i+1}] \longleftarrow \ldots$$

*is almost-exact. Then the induced sequence of abelian varieties*

$$\ldots \longrightarrow J_{H_{i-1}} \xrightarrow{\phi_{i-1}} J_{H_i} \xrightarrow{\phi_i} J_{H_{i+1}} \longrightarrow \ldots$$

*is almost-exact. Furthermore, we have*

$$\frac{\#\Phi(J_{H_i})(\mathbb{C})}{\#\Phi(\mathbb{Z}[G/H_i])} \leq \#\eta(\mathbb{Z}[G/H_{i-1}], \operatorname{im}\sigma_{i-1}) \cdot \#\ker \pi_{H_{i+1}*}$$
$$\cdot |H_i| \cdot \#\operatorname{coker} \pi_{H_i}^* \cdot \#\frac{\ker \pi_{H_i}^*}{\ker \pi_{H_i}^* \cap \operatorname{im} \phi_{i-1}}.$$

*Proof.* The almost-exact sequence

$$\ldots \longleftarrow \mathbb{Z}[G/H_{i-1}] \xleftarrow{\sigma_{i-1}} \mathbb{Z}[G/H_i] \xleftarrow{\sigma_i} \mathbb{Z}[G/H_{i+1}] \longleftarrow \ldots$$

induces by Lemma 3.4 an almost-exact sequence

$$\ldots \longrightarrow J(\mathbb{C})^{H_{i-1}} \xrightarrow{\phi'_{i-1}} J(\mathbb{C})^{H_i} \xrightarrow{\phi'_i} J(\mathbb{C})^{H_{i+1}} \longrightarrow \ldots$$

where $\#\Phi(J(\mathbb{C})^{H_i}) \leq \#\Phi(\mathbb{Z}[G/H_i]) \cdot \#\eta(\mathbb{Z}[G/H_{i-1}], \operatorname{im}\sigma_{i-1})$.

By Lemma 3.5, we have a commutative diagram

$$\begin{array}{ccccccc}
J_{H_{i-1}}(\mathbb{C}) & \xrightarrow{\phi_{i-1}} & J_{H_i}(\mathbb{C}) & \xrightarrow{=} & J_{H_i}(\mathbb{C}) & \xrightarrow{\phi_i} & J_{H_{i+1}}(\mathbb{C}) \\
\pi_{H_{i-1}}^* \downarrow & & \pi_{H_i *} \uparrow & & \pi_{H_i}^* \downarrow & & \pi_{H_{i+1}*} \uparrow \\
J(\mathbb{C})^{H_{i-1}} & \xrightarrow{\phi'_{i-1}} & J(\mathbb{C})^{H_i} & \xrightarrow{\pi_{H_i}^* \circ \pi_{H_i *}} & J(\mathbb{C})^{H_i} & \xrightarrow{\phi'_i} & J(\mathbb{C})^{H_{i+1}}.
\end{array}$$

It is easily seen that $\phi_i \circ \phi_{i-1} = 0$ using the fact that $\phi'_i \circ \phi'_{i-1} = 0$ and $\pi_{H_i}^* \circ \pi_{H_i *}$ is multiplication by $\deg(\pi_{H_i}) = |H_i|$ when restricted to the subgroup variety $J^{H_i}(\mathbb{C}) = J(\mathbb{C})^{H_i}$, so the image of $\phi'_{i-1}$ is sent into itself under this morphism.

We next show the relation between $\#\Phi(J_{H_i}(\mathbb{C}))$ and $\#\Phi(J(\mathbb{C})^{H_i})$. To begin with, note that

$$\frac{\pi_{H_i}^*(\ker \phi_i)}{\pi_{H_i}^*(\operatorname{im} \phi_{i-1})} \cong \frac{\ker \phi_i}{\ker \phi_i \cap (\operatorname{im} \phi_{i-1} + \ker \pi_{H_i}^*)}$$

so that

$$\#\Phi(J_{H_i}(\mathbb{C})) = \#\frac{\ker \phi_i}{\operatorname{im} \phi_{i-1}} = \#\frac{\pi^*_{H_i}(\ker \phi_i)}{\pi^*_{H_i}(\operatorname{im} \phi_{i-1})} \cdot \#\frac{\ker \phi_i \cap (\operatorname{im} \phi_{i-1} + \ker \pi^*_{H_i})}{\operatorname{im} \phi_{i-1}}$$

$$= \#\frac{\pi^*_{H_i}(\ker \phi_i)}{\pi^*_{H_i}(\operatorname{im} \phi_{i-1})} \cdot \#\frac{\operatorname{im} \phi_{i-1} + \ker \pi^*_{H_i}}{\operatorname{im} \phi_{i-1}}$$

$$= \#\frac{\pi^*_{H_i}(\ker \phi_i)}{\pi^*_{H_i}(\operatorname{im} \phi_{i-1})} \cdot \#\frac{\ker \pi^*_{H_i}}{\ker \pi^*_{H_i} \cap \operatorname{im} \phi_{i-1}}$$

where the second last equality follows from the fact that $A \cap (B+C) = A \cap B + A \cap C$ if $B \subset A$, $\operatorname{im} \phi_{i-1} \subset \ker \phi_i$, and $\ker \pi^*_{H_i} \subset \ker \phi_i$.

By the commutative diagram, we see that

$$\pi^*_{H_i}(\ker \phi_i) = \ker \pi_{H^*_{i+1}} \circ \phi'_i$$
$$\pi^*_{H_i}(\operatorname{im} \phi_{i-1}) = \operatorname{im}[|H_i|] \circ \phi'_{i-1} \circ \pi^*_{H_{i-1}}.$$

Now, we have

$$\#\frac{\ker \pi_{H^*_{i+1}} \circ \phi'_i}{\ker \phi'_i} = \#\ker \pi_{H^*_{i+1}}$$

$$\#\frac{\operatorname{im}[|H_i|] \circ \phi'_{i-1}}{\operatorname{im}[|H_i|] \circ \phi'_{i-1} \circ \pi^*_{H_{i-1}}} = \#\operatorname{coker} \pi^*_{H_{i-1}}.$$

This yields

$$\#\frac{\pi^*_{H_i}(\ker \phi_i)}{\pi^*_{H_i}(\operatorname{im} \phi_{i-1})} = \#\frac{\ker \phi'_i}{\operatorname{im}[|H_i|] \circ \phi'_{i-1}} \cdot \#\ker \pi_{H^*_{i+1}} \cdot \#\operatorname{coker} \pi^*_{H_{i-1}}.$$

Since

$$\#\frac{\ker \phi'_i}{\operatorname{im}[|H_i|] \circ \phi'_{i-1}} = \frac{\ker \phi'_i}{\operatorname{im} \phi'_{i-1}} \cdot |H_i|$$

we obtain finally that

$$\frac{\#\Phi(J_{H_i}(\mathbb{C}))}{\#\Phi(J(\mathbb{C})^{H_i})} \leq \#\ker \pi_{H^*_{i+1}} \cdot |H_i| \cdot \#\operatorname{coker} \pi^*_{H_{i-1}} \cdot \#\frac{\ker \pi^*_{H_i}}{\ker \pi^*_{H_i} \cap \operatorname{im} \phi_{i-1}}.$$

The result then follows from Lemma 3.6 and 3.2. □

## 4. Double cosets

In this section, an explicit description of the $\mathbb{Z}[G]$-modules homomorphisms between two induced representations $\mathbb{Z}[G/H]$ and $\mathbb{Z}[H/K]$ is given in terms of double cosets. This will be used to describe $\mathbb{Z}[G]$-modules homomorphisms which induce the homomorphisms $(\phi_{GB}, \phi_{BN}, \phi_{NN'})$ via Lemma 3.2.

Let $H, K$ be subgroups of $G$. Suppose we have a $\mathbb{Z}[G]$-module homomorphism $\phi : \mathbb{Z}[G/H] \to \mathbb{Z}[G/K]$. Then $\phi$ is determined by its value on the coset $H$ since $\phi(gH) = g\phi(H)$. Moreover, for all $h \in H$ we have $h\phi(H) = \phi(hH) = \phi(H)$ so the element $\phi(H) \in \mathbb{Z}[G/K]$ is invariant under multiplication on the left by $h$.

Conversely, suppose we are given an element $\alpha$ of $\mathbb{Z}[G/K]$ which is invariant under multiplication on the left by elements in $H$. Then one can define a $\mathbb{Z}[G]$-module homomorphism $\phi : \mathbb{Z}[G/H] \to \mathbb{Z}[G/K]$ such that $\phi(H) = \alpha$. Indeed, define $\phi(gH) = g\alpha$. This is well-defined on right $H$-cosets because $\alpha$ is invariant on the left by multiplication by elements in $H$. The desired map is then obtained by extending this $\mathbb{Z}$-linearly to all of $\mathbb{Z}[G/H]$. By its very definition, it is a $\mathbb{Z}[G]$-module homomorphism.

**Lemma 4.1.**
$$HgK = \bigcup_{\alpha \in H/H_g} \alpha gK$$

where $H_g = H \cap gKg^{-1}$ and the union is disjoint. We call $[H : H_g]$ the degree of $HgK$.

*Proof.* Since
$$H = \bigcup_{\alpha \in H/H_g} \alpha H_g$$

we have
$$HgK = \bigcup_{\alpha \in H/H_g} \alpha H_g gK$$
$$= \bigcup_{\alpha \in H/H_g} \alpha gK.$$

If $\alpha gK = \alpha' gK$, then $\alpha g \in \alpha' gK$ and hence $\alpha \in \alpha' H_g$. □

There is a distinguished class of $\mathbb{Z}[G]$-module homomorphisms from $\mathbb{Z}[G/H]$ to $\mathbb{Z}[G/K]$ : for each double coset $HgK$ where $g \in G$, define $\phi(HgK)(H) = HgK$, where the double coset $HgK$ is considered as an element of $\mathbb{Z}[G/K]$ by decomposing it into right $K$-cosets. This yields a $\mathbb{Z}[G]$-module homomorphism $\phi(HgK) : \mathbb{Z}[G/H] \to \mathbb{Z}[G/K]$ as above.

**Lemma 4.2.** *Let $H, K$ be subgroups of $G$. Then as $\mathbb{Z}$-modules*
$$\Theta : \mathbb{Z}[H\backslash G/K] \cong \mathrm{Hom}_{\mathbb{Z}[G]}(\mathbb{Z}[G/H], \mathbb{Z}[G/K]).$$

*Proof.* Let $\Omega$ be a complete set of inequivalent representatives for $H\backslash G/K$. It is clear that $\phi$ is injective since $\Theta(\sum_{g\in\Omega} \alpha_g HgK) = 0$ means that $\sum_{g\in\Omega} \alpha_g HgK = 0$ as an element of $\mathbb{Z}[G/K]$. This occurs if and only if $\alpha_g = 0$ for all $g \in \Omega$.

A $\mathbb{Z}[G]$-module homomorphism $\Theta : \mathbb{Z}[G/H] \to \mathbb{Z}[G/K]$ is determined by its value on the coset $H$. Since $\Theta(H)$ is an element in $\mathbb{Z}[G/K]$ which is invariant under multiplication by elements in $H$, we can write $\Theta(H) = \sum_{g\in\Omega} \alpha_g HgK$. We then see that $\Theta = \sum_{g\in\Omega} \alpha_g \Theta(HgK)$. This shows $\Theta$ is surjective. □

We call elements of $\mathbb{Z}[H\backslash G/K]$ operators and by convention denote their action from the left. This causes the slight annoyance that $HK \times KH = KH \circ HK$. We will also omit $\Theta$ when considering an operator as an element of $\mathrm{Hom}_{\mathbb{Z}[G]}(\mathbb{Z}[G/H], \mathbb{Z}[G/K])$.

In the case $H = K$, $\mathrm{End}_{\mathbb{Z}[G]}(\mathbb{Z}[G/H])$ is not just a $\mathbb{Z}$-module but also a ring with unit. Thus $\Theta$ gives a ring structure on $\mathbb{Z}[H\backslash G/H]$. This ring structure can be described explictly and developed purely in terms of $\mathbb{Z}[H\backslash G/H]$ (c.f. [29] and below). Finally, we remark that $\mathbb{Z}[H\backslash G/H]$ is also nothing other than the Hecke algebra $\mathcal{H}(G, H)$ with trivial character which arises in representation theory [6], though it is used in a different context there.

RELATIONS BETWEEN JACOBIANS 9

More generally, the $\mathbb{Z}$-linear multiplication

$$\mathrm{Hom}_{\mathbb{Z}[G]}(\mathbb{Z}[G/H], \mathbb{Z}[G/K]) \times \mathrm{Hom}_{\mathbb{Z}[G]}(\mathbb{Z}[G/K], \mathbb{Z}[G/M])$$
$$\to \mathrm{Hom}_{\mathbb{Z}[G]}(\mathbb{Z}[G/H], \mathbb{Z}[G/M])$$
$$(f, g) \mapsto g \circ f$$

induces via $\Theta$ a $\mathbb{Z}$-linear multiplication

$$\mathbb{Z}[H\backslash G/K] \times \mathbb{Z}[K\backslash G/M] \to \mathbb{Z}[H\backslash G/M].$$

It is easy to describe this multiplication via $\Theta$. For example, suppose

$$HaK = \cup_{\alpha \in \Omega_a} \alpha aK$$
$$KbM = \cup_{\beta \in \Omega_b} \beta bM.$$

Then

$$HaK \times KbM = \sum_{\alpha \in \Omega_a} \sum_{\beta \in \Omega_b} \alpha a \beta b M = \sum_{g \in \Omega} \gamma_g HgM.$$

here $\Omega$ is a complete set of inequivalent representatives for $H\backslash G/M$. Note that $\gamma_g$ is easily calculated to be

$$\gamma_g = \#\{\alpha a \beta b \in HgM\} / \deg(HgM)$$
$$= \#\{\alpha a \beta b \in gM\}$$

(compare with the formula given in [29]). In our context, all groups are finite so this can be alternatively described as

$$HaK \times KbM = \frac{\deg(KbM)}{|K|} \sum_{k \in K} \frac{\deg(HaK)}{\deg(HakbM)} HakbM.$$

## 5. Decomposition into irreducibles

In this section, we briefly review the representation theory of $G$ and decompose the induced representations $\mathbb{C}[G/G]$, $\mathbb{C}[G/B]$, $\mathbb{C}[G/N]$, $\mathbb{C}[G/N']$. This material is standard and can be found in [12] for instance.

The set of conjugacy classes $C(G)$ of $G$ can be described by breaking it into four types

$$C(G) = \left\{ [h_x = \begin{pmatrix} x & 0 \\ 0 & x \end{pmatrix}] \mid x \in \mathbb{F}_p^\times \right\} \bigcup$$
$$\left\{ [b_x = \begin{pmatrix} x & 1 \\ 0 & x \end{pmatrix}] \mid x \in \mathbb{F}_p^\times \right\} \bigcup$$
$$\left\{ [\kappa_{x,y} = \begin{pmatrix} x & 0 \\ 0 & y \end{pmatrix}] \mid (x,y) \in \mathbb{F}_p^\times \times \mathbb{F}_p^\times - \Delta \text{ and } (x,y) \sim (y,x) \right\} \bigcup$$
$$\left\{ [\gamma_{x,y} = \begin{pmatrix} x & \lambda y \\ y & x \end{pmatrix}] \mid (x,y) \neq (0,0) \in \mathbb{F}_p \times \mathbb{F}_p \text{ and } (x,y) \sim (x,-y) \right\},$$

where we use the notation $[g]$ to denote the conjugacy class of $g \in G$, and $\sim$ the relation of conjugation. The elements $h_x$, $b_x$, $\kappa_{x,y}$, $\gamma_{x,y}$ give distinct conjugacy classes except for the identifications $\sim$ given above.

Let $\alpha$, $\beta$, and $\phi$ be one-dimensional representations of $\mathbb{F}_p^\times$, $\mathbb{F}_p^\times$, and $\mathbb{F}_{p^2}^\times$, respectively. There are four types of irreducible representations of $G$. They are denoted $U_\alpha$, $V_\alpha$, $W_{\alpha,\beta}$, and $X_\phi$, with the restrictions $\alpha \neq \beta$ and $\phi^p \neq \phi$. These representations are all distinct except for the isomorphisms $W_{\alpha,\beta} \cong W_{\beta,\alpha}$ and $X_\phi \cong X_{\phi^p}$. Their values on each conjugacy class is given in the following table.

TABLE 1. The character table of $G = \mathrm{GL}_2(\mathbb{F}_p)$

|  | $[h_x]$ | $[b_x]$ | $[\kappa_{x,y}]$ | $[\gamma_{x,y}] = [\gamma]$ |
|---|---|---|---|---|
| $U_\alpha$ | $\alpha(x^2)$ | $\alpha(x^2)$ | $\alpha(xy)$ | $\alpha(\gamma^{p+1})$ |
| $V_\alpha$ | $p\alpha(x^2)$ | $0$ | $\alpha(xy)$ | $-\alpha(\gamma^{p+1})$ |
| $W_{\alpha,\beta}$ | $(p+1)\alpha(x)\beta(x)$ | $\alpha(x)\beta(x)$ | $\alpha(x)\beta(y)+\alpha(y)\beta(x)$ | $0$ |
| $X_\phi$ | $(p-1)\phi(x)$ | $-\phi(x)$ | $0$ | $-(\phi(\gamma)+\phi(\gamma^p))$ |

A routine calculation using the above character table gives the following propositions.

**Proposition 5.1.**
$$\mathbb{C}[G/G] \cong U_1$$

*Proof.* This is clear. □

**Proposition 5.2.**
$$\mathbb{C}[G/B] \cong U_1 \oplus V_1$$

**Proposition 5.3.**

$$\mathbb{C}[G/N'] \cong U_1 \oplus \sum_{p \equiv 1 \pmod{4}, \alpha = \left(\frac{\cdot}{p}\right)} V_\alpha \oplus \sum_{\alpha^{\frac{p-1}{2}}=1, \alpha^2 \neq 1} W_{\alpha,\alpha^{-1}} \oplus \sum_{\phi^{p+1}=1, \phi^{\frac{p+1}{2}} \neq 1, \phi^{p-1} \neq 1} X_\phi.$$

**Proposition 5.4.**

$$\mathbb{C}[G/N] \cong U_1 \oplus \sum_{p \equiv 1 \pmod{4}, \alpha = \left(\frac{\cdot}{p}\right)} V_\alpha \oplus \sum_{\alpha^{\frac{p-1}{2}}=1, \alpha^2 \neq 1} W_{\alpha,\alpha^{-1}} \oplus \sum_{\phi^{p+1}=1, \phi^{\frac{p+1}{2}} \neq 1, \phi^{p-1} \neq 1} X_\phi \oplus V_1$$

*Proof.* This follows from the previous propositions as
$$\mathbb{C}[G/N'] \oplus \mathbb{C}[G/B] \cong \mathbb{C}[G/N] \oplus \mathbb{C}[G/G]$$
by Theorem 4. □

For later reference, we note that

(9) $$\phi^{p+1} = 1 \iff \phi\mid_{\mathbb{F}_p^\times} = 1$$

(10) $$\alpha^{\frac{p-1}{2}} = 1 \iff \alpha(-1) = 1.$$

## 6. Exactness at $\mathbb{C}[G/G]$ and $\mathbb{C}[G/B]$

In this section, we exhibit choices for the $\mathbb{Z}[G]$-modules homomorphisms occuring in sequence 3 which via Lemma 3.2 induce the homomorphisms of jacobians in Theorem 2. An analysis of the almost-exactness of this sequence of $\mathbb{Z}[G]$-modules is performed. In particular, we reduce the property of almost-exactness to the calculation of certain eigenvalues. The first two terms in the sequence of $\mathbb{Z}[G]$-modules is easily shown to be exact.

Before starting, we state some degree of operators for later reference.

**Lemma 6.1.**
$$\deg(NN') = (p-1)/2$$
$$\deg(N'N) = (p+1)/2$$

**Lemma 6.2.**
$$\deg(NB) = 2$$
$$\deg(BN) = p$$

**Lemma 6.3.**
$$\deg(NG) = 1$$
$$\deg(GN) = p(p+1)/2$$

**Lemma 6.4.**
$$\deg(BG) = 1$$
$$\deg(GB) = p+1$$

The following lemma exhibits the natural operators coming from representation theory which induce the homomorphisms of jacobians in Theorem 2.

**Lemma 6.5.** *The $\mathbb{Z}[G]$-module homomorphisms*
$$(\sigma_{BG}, \sigma_{NB}, \sigma_{N'N}) = (BG, |G|\,(1 - \mathrm{pr}_G)NB, N'N)$$
*induce via Lemma 3.2 the homomorphisms $(\phi_{GB}, \phi_{BN}, \phi_{NN'})$ in Theorem 2.*

*Proof.* This is easily checked by going through Lemma 3.2. Indeed, we have the following verifications.

$$\sigma_{BG} = \pi_{B*} \circ BG(B) \circ \pi_G^*$$
$$= \pi_{B*} \circ \pi_G^*$$

$$\sigma_{NB} = \pi_{N*} \circ (|G|\,(1 - \mathrm{pr}_G)NB)(N) \circ \pi_B^*$$
$$= \pi_{N*} \circ |G|\,(1 - \mathrm{pr}_G)(1 + \omega) \circ \pi_B^*$$

$$\sigma_{N'N} = \pi_{N'*} \circ N'N(N') \circ \pi_N^*$$
$$= \left[\frac{p+1}{2}\right] \circ \pi_{N'*} \circ \pi_N^*$$

where $N'N = \cup_{\alpha \in \Omega} \alpha N$ is a disjoint union, and we have used the fact that $\deg(N'N) = \frac{p+1}{2}$ shown in Lemma 6.1. □

The above lemma together with Proposition 3.7, implies that Theorem 2 would follow if one can show that the sequence of $\mathbb{Z}[G]$-modules

(11)
$$0 \longleftarrow \mathbb{Z}[G/G] \xleftarrow{\sigma_{BG}} \mathbb{Z}[G/B] \xleftarrow{\sigma_{NB}} \mathbb{Z}[G/N] \xleftarrow{\sigma_{N'N}} \mathbb{Z}[G/N'] \longleftarrow 0$$

is almost-exact. To do this, we essentially dualise the problem. The basis for this is given in next two lemmas.

**Lemma 6.6.** *Suppose we have a sequence*

$$\ldots \longleftarrow \mathbb{Z}[G/H_{i-1}] \xleftarrow{\sigma_{i-1}} \mathbb{Z}[G/H_i] \xleftarrow{\sigma_i} \mathbb{Z}[G/H_{i+1}] \longleftarrow \ldots$$

*together with $\mathbb{Z}[G]$-module homomorphisms $\tau_i : \mathbb{Z}[G/H_i] \to \mathbb{Z}[G/H_{i+1}]$, where $\ker \sigma_{i-1} \supset \operatorname{im} \sigma_i$. Then the above sequence is almost-exact if $\ker \tau_{i-1} \circ \sigma_{i-1} \supset \operatorname{im} \sigma_i \circ \tau_i$ with finite quotient.*

*Proof.* Simply note that

$$\ker \tau_{i-1} \circ \sigma_{i-1} \supset \ker \sigma_{i-1} \supset \operatorname{im} \sigma_i \supset \operatorname{im} \sigma_i \circ \tau_i.$$

□

Let $V$ be $\mathbb{C}[G]$-module. We say $V$ has multiplicity one if every irreducible component occurs with multiplicity at most one. If this is the case, then a $\mathbb{C}[G]$-module homomorphism $\sigma : V \to V$ is given by a scalar on each irreducible component of $V$ by Schur's lemma. If the irreducible component has character $\chi$, then this eigenvalue is given by

$$\lambda_\chi(\sigma) = \frac{1}{\chi(1)} \operatorname{tr}(\operatorname{pr}_\chi \circ \sigma) = \frac{1}{\chi(1)} \operatorname{tr}(\sigma \circ \operatorname{pr}_\chi)$$

where $\operatorname{pr}_\chi = \frac{\chi(1)}{|G|} \sum_{g \in G} \overline{\chi(g)} g$ is the projector to the $\chi$-component.

**Lemma 6.7.** *Suppose we have a sequence of $\mathbb{Z}[G]$-modules*

$$M \xleftarrow{\epsilon} M \xleftarrow{\delta} M$$

*where $\ker \epsilon \supset \operatorname{im} \delta$ and $M \otimes \mathbb{C}$ has multiplicity one. Then the quotient $\Phi = \ker \epsilon / \operatorname{im} \delta$ is finite if and only if $\lambda_\chi(\epsilon \otimes \mathbb{C}) = 0$ implies $\lambda_\chi(\delta \otimes \mathbb{C}) \neq 0$ for every character $\chi$ of an irreducible component of $M$. If this is the case then*

$$\#\Phi(M) = \prod_{(\chi,\chi_M)=1, \lambda_\chi(\delta \otimes \mathbb{C}) \neq 0} \lambda_\chi(\delta \otimes \mathbb{C})^{\chi(1)}$$

*Proof.* To show the first part of the lemma, consider the sequence

$$L \xleftarrow{\epsilon \otimes \mathbb{C}} L \xleftarrow{\delta \otimes \mathbb{C}} L$$

where $L = M \otimes \mathbb{C}$. Decompose $L = \oplus_\chi L_\chi$, where $L_\chi$ is the $\chi$-component of $L$, that is the direct sum of the irreducible components of $L$ with character $\chi$. By the multiplicity one assumption, $L_\chi$ is irreducible. We use the convention that $\chi$ runs through all irreducible characters of $G$ so that $L_\chi$ can be zero for some $\chi$.

By Schur's lemma, the $\mathbb{C}[G]$-module homomorphism $\delta_i \otimes \mathbb{C} : L \to L$ is multiplication by a scalar $\lambda_\chi(\delta \otimes \mathbb{C})$ on each irreducible component $L_\chi$. We have that

$$\ker \epsilon \otimes \mathbb{C} = \prod_{(\chi,\chi_M)=1, \lambda_\chi(\epsilon \otimes \mathbb{C})=0} L_\chi$$

$$\operatorname{im} \delta \otimes \mathbb{C} = \prod_{(\chi,\chi_M)=1, \lambda_\chi(\delta \otimes \mathbb{C}) \neq 0} L_\chi.$$

By hypothesis, $\ker \epsilon \otimes \mathbb{C} \supset \operatorname{im} \delta \otimes \mathbb{C}$. However, if $\lambda_\chi(\epsilon \otimes \mathbb{C}) = 0$ implies $\lambda_\chi(\delta \otimes \mathbb{C}) \neq 0$, then the above descriptions of $\ker \epsilon \otimes \mathbb{C}$ and $\operatorname{im} \delta \otimes \mathbb{C}$ show that they are in fact equal.

Consider the given sequence again. Replacing the middle term by $\ker \delta_{i-1}$, we have a sequence

$$M \xleftarrow{\epsilon} \ker \epsilon \xleftarrow{\delta} M$$

where $\Phi(\ker \epsilon) = \Phi(M)$. However, $\Phi(\ker \epsilon)$ is simply the cokernel of the map of $\mathbb{Z}[G]$-module homomorphism $\delta$. The order of this cokernel is given by the determinant of the map $\delta \otimes \mathbb{C}$ (considered as a map to $\ker \epsilon \otimes \mathbb{C} = \operatorname{im} \delta \otimes \mathbb{C}$), which is easily seen to be $\prod_{(\chi, \chi_M)=1, \lambda_\chi(\delta \otimes \mathbb{C}) \neq 0} \lambda_\chi(\delta \otimes \mathbb{C})^{\chi(1)}$. □

**Corollary 6.8.** *Suppose we have a sequence*

$$\ldots \longleftarrow \mathbb{Z}[G/H_{i-1}] \xleftarrow{\sigma_{i-1}} \mathbb{Z}[G/H_i] \xleftarrow{\sigma_i} \mathbb{Z}[G/H_{i+1}] \longleftarrow \ldots$$

*together with $\mathbb{Z}[G]$-module homomorphisms $\tau_i : \mathbb{Z}[G/H_i] \to \mathbb{Z}[G/H_{i+1}]$, where $\ker \sigma_{i-1} \supset \operatorname{im} \sigma_i$, and each $\mathbb{C}[G/H_i]$ has multiplicity one. Put $\delta_i = \sigma_i \circ \tau_i$ and $\epsilon_i = \tau_{i+1} \circ \sigma_{i+1}$. If $\lambda_\chi(\epsilon_i \otimes \mathbb{C}) = 0$ implies $\lambda_\chi(\delta_i \otimes \mathbb{C}) \neq 0$ for every character $\chi$ of an irreducible component of $\mathbb{C}[G/H_i]$, then the original sequence is almost-exact at $\mathbb{Z}[G/H_i]$, and $\#\Phi(\mathbb{Z}[G/H_i]) = \prod_{\lambda_\chi(\delta_i \otimes \mathbb{C}) \neq 0} \lambda_\chi(\delta_i \otimes \mathbb{C})^{\chi(1)}$.*

We wish to apply the above corollary to sequence 11 by taking the $\mathbb{Z}[G]$-module homomorphisms $(\tau_{GB}, \tau_{BN}, \tau_{NN'})$ to be $(GB, |G|(1 - \operatorname{pr}_G)BN, NN')$. There are several conditions to check.

**Lemma 6.9.** *We have that $\sigma_{N'B} \circ \sigma_{NN'} = 0$ and $\sigma_{BG} \circ \sigma_{N'B} = 0$.*

*Proof.* Using the fact that $N'B = G$, we see that $|G|\operatorname{pr}_G \circ N'B = |G| N'B$, which shows that

$$\sigma_{NB} \circ \sigma_{NN'} = |G|(1 - \operatorname{pr}_G) \circ N'N \times NB$$
$$= 0.$$

Similarly, using the fact that $NG = G$, we see that $|G|\operatorname{pr}_G \circ = |G| NG$, which shows that

$$\sigma_{BG} \circ \sigma_{NB} = BG \circ |G|(1 - \operatorname{pr}_G)NB$$
$$= |G|(1 - \operatorname{pr}_G) \circ NB \times BG$$
$$= 0.$$

□

**Lemma 6.10.** *The $\mathbb{C}[G]$-modules $\mathbb{C}[G/G]$, $\mathbb{C}[G/B]$, $\mathbb{C}[G/N]$, $\mathbb{C}[G/N']$ are of multiplicity one.*

*Proof.* This lemma follows from the calculations in Section 5. □

To check the remaining hypotheses of Corollary 6.8, one needs to compare the eigenvalues of

(12) $\quad \epsilon_G = 0$

(13) $\quad \epsilon_B = \tau_{GB} \circ \sigma_{BG} = GB \circ BG = BG \times GB$

(14) $\quad \epsilon_N = \tau_{BN} \circ \sigma_{NB} = |G|(1 - \operatorname{pr}_G)BN \circ |G|(1 - \operatorname{pr}_G)NB$

(15) $\quad \phantom{\epsilon_N} = |G|^2 (1 - \operatorname{pr}_G)NB \times BN$

(16) $\quad \epsilon_{N'} = \tau_{NN'} \circ \sigma_{N'N} = NN' \circ N'N = N'N \times NN'$

with the eigenvalues of

(17) $$\delta_G = \sigma_{BG} \circ \tau_{GB} = BG \circ GB = GB \times BG$$

(18) $$\delta_B = \sigma_{NB} \circ \tau_{BN} = |G|\,(1 - \mathrm{pr}_G)NB \circ |G|\,(1 - \mathrm{pr}_G)BN$$

(19) $$= |G|^2\,(1 - \mathrm{pr}_G)BN \times NB$$

(20) $$\delta_N = \sigma_{N'N} \circ \tau_{NN'} = N'N \circ NN' = NN' \times N'N$$

(21) $$\delta_{N'} = 0,$$

respectively.

The eigenvalue $\lambda_\chi(|G|^2\,(1 - \mathrm{pr}_G))$ is simply 0 if $\chi$ is the trivial character and $|G|^2$ otherwise. Hence, we are led to calculating eigenvalues of operators of the form $HK \times KH$, where $H, K$ are subgroups of $G$. The following lemma gives an expression for such an eigenvalue.

**Lemma 6.11.** *Let $\chi$ be an irreducible character of $G$ and $H, K$ subgroups of $G$. Then the trace of $HK \times KH$ on the $\chi$-component of $\mathbb{C}[G/H]$ is given by*

$$\mathrm{tr}_\chi(HK \times KH) = \chi(1)\frac{\deg(HK)}{|H|}\frac{\deg(KH)}{|K|}\sum_{k \in K}\sum_{h \in H}\overline{\chi(kh)}$$

*Proof.* Note that

$$\mathrm{tr}_\chi(HK \times KH) = \mathrm{tr}(\mathrm{pr}_\chi \circ HK \times KH).$$

Choose a set of inequivalent representatives $g_1, \ldots, g_n$ for $G/H$, where $n = |G/H|$. Then $g_1 H, \ldots, g_n H$ forms a $\mathbb{C}$-basis for $\mathbb{C}[G/H]$. To calculate the trace of the map $\mathrm{pr}_\chi \circ HK \times KH$, it suffices to compute for each $i$, the coefficient $\alpha_i$ of $g_i H$ in $(\mathrm{pr}_\chi \circ HK \times KH)(g_i H)$. The trace is then given by $\sum_{i=1}^n \alpha_i$.

To begin with, we have

$$\mathrm{pr}_\chi \circ HK \times KH(g_i H) = \left(\frac{\chi(1)}{|G|}\sum_{g \in G}\overline{\chi(g)}g\right) \cdot \left(\frac{\deg(HK)}{|H|}\frac{\deg(KH)}{|K|}\sum_{h \in H}\sum_{k \in K}g_i hk H\right)$$

$$= \frac{\chi(1)}{|G|} \cdot \frac{\deg(HK)}{|H|} \cdot \frac{\deg(KH)}{|K|}\sum_{h \in H}\sum_{k \in K}\sum_{g \in G}\overline{\chi(g)}g_i hkgH$$

where the last equality follows from the fact that $\mathrm{pr}_\chi$ is a $\mathbb{C}[G]$-module homomorphism. The coefficient of $g_i H$ in the above element of $\mathbb{C}[G/H]$ is then given by

$$\alpha_i = \frac{\chi(1)}{|G|} \cdot \frac{\deg(HK)}{|H|} \cdot \frac{\deg(KH)}{|K|} \cdot |H|\sum_{k \in K}\sum_{h \in H}\overline{\chi(kh)}.$$

Thus,

$$\mathrm{tr}_\chi(HK \times KH) = \chi(1)\frac{\deg(HK)}{|H|} \cdot \frac{\deg(KH)}{|K|}\sum_{k \in K}\sum_{h \in H}\overline{\chi(kh)}$$

as desired. $\square$

Furthermore, we have that

**Lemma 6.12.** *Let $\chi$ be an irreducible character of $G$ and $H, K$ subgroups of $G$. Then*

$$\mathrm{tr}_\chi(HK \times KH) = \mathrm{tr}_\chi(KH \times HK).$$

*Proof.* Since the elements $hk$ and $kh$ are conjugate, the result follows from comparing the expressions for $\text{tr}_\chi(HK \times KH)$ and $\text{tr}_\chi(KH \times HK)$ in Lemma 6.11. □

We now compute some eigenvalues.

**Lemma 6.13.** *Let $\chi$ be an irreducible character of $G$. Then*

$$\lambda_\chi(BN \times NB) = \lambda_\chi(NB \times BN) = \begin{cases} 2p & \text{if } \chi = U_1 \\ p^2 + p - 1 & \text{if } \chi = V_1 \\ 0 & \text{otherwise} \end{cases}.$$

*Proof.* Note that $\text{tr}_\chi(BN \times NB) = 0$ if $\chi$ is not the character of $U_1$ nor $V_1$ since these are the only two irreducibles which occur in $\mathbb{C}[G/B]$ (see Lemma 5.2).

From Lemma 6.11, we see that

$$\text{tr}_\chi(BN \times NB) = \chi(1) \frac{\deg(BN)}{|B|} \frac{\deg(NB)}{|N|} \sum_{n \in N} \sum_{b \in B} \overline{\chi(bn)}.$$

For $\chi = U_1$, this implies that

$$\text{tr}_\chi(BN \times NB) = \chi(1) \frac{\deg(BN)}{|B|} \frac{\deg(NB)}{|N|} |B| \, |N| = 2p^2.$$

For $\chi = V_1$, we simplify a bit further and see that

$$\text{tr}_\chi(BN \times NB) = \chi(1) \frac{\deg(BN)}{|B|} \frac{\deg(NB)}{|N|} \sum_{c \in C} \sum_{b \in B} (\overline{\chi(bc)} + \overline{\chi(bc\omega)})$$

$$= \chi(1) \frac{\deg(BN)}{|B|} \frac{\deg(NB)}{|N|} |C| \sum_{b \in B} (\overline{\chi(b)} + \overline{\chi(b\omega)}).$$

Now, note that $\sum_{b \in B} \overline{\chi(b)} = |G| \, (\chi, 1_B) = |G|$. On the other hand,

$$\sum_{b \in B} \overline{\chi(b\omega)} = \sum_{[g] \in C(G)} \overline{\chi([g])} c([g])$$

where $C(G)$ is the set of all conjugacy classes, $[g]$ is the conjugacy class of $g \in G$, and $c([g])$ denotes the number of elements in $B\omega$ which lie in $[g]$. The quantity $c([g])$ is easily calculated to be $p-1$ if $[g]$ is non-scalar and 0 otherwise by counting the number of elements in $B\omega$ with fixed trace and determinant. From the values of $\chi$ on each conjugacy class (see section 5), we see that

$$\sum_{[g] \in C(G)} \overline{\chi([g])} c([g]) = ((p-1)(p-2)/2 - (p^2-p)/2)(p-1) = -(p-1)^2.$$

Thus, we have that

$$\text{tr}_\chi(BN \times NB) = \chi(1) \frac{\deg(BN)}{|B|} \frac{\deg(NB)}{|N|} |C| \, [(p^2-p)(p^2-1) - (p-1)^2]$$

$$= p \frac{2}{p(p-1)^2} \frac{p}{2(p-1)^2} (p-1)^2 (p-1)^2 (p^2+p-1)$$

$$= p(p^2 + p - 1)$$

and hence

$$\lambda_\chi(BN \times NB) = p^2 + p - 1$$

□

**Lemma 6.14.** *Let $\chi$ be an irreducible character of $G$. Then*

$$\lambda_\chi(GB \times BG) = \lambda_\chi(BG \times GB) = \begin{cases} p+1 & \text{if } \chi = U_1 \\ 0 & \text{otherwise} \end{cases}.$$

*Proof.* Note that $\text{tr}_\chi(GB \times BG) = 0$ if $\chi$ is not the character of $U_1$ since $\mathbb{C}[G/G]$ is the trivial representation. (see Lemma 5.1).

For $\chi = U_1$, we see from Lemma 6.11 that

$$\text{tr}_\chi(GB \times BG) = \chi(1) \frac{\deg(GB)}{|G|} \frac{\deg(BG)}{|B|} \sum_{b \in B} \sum_{g \in G} \overline{\chi(bg)}$$
$$= \deg(GB) \deg(BG)$$
$$= p + 1.$$

□

**Proposition 6.15.** *Sequence 11 is exact at $\mathbb{Z}[G/G]$ and $\mathbb{Z}[G/B]$.*

*Proof.* Since $\mathbb{C}[G/G]$ is the trivial representation, we only need to check the hypotheses of Corollary 6.8 for $\chi = U_1$. Indeed, $\lambda_\chi(\epsilon_G) = 0$ and $\lambda(\delta_G) = p + 1 \neq 0$ by Lemma 6.14. Thus, by the corollary, sequence 11 is exact at $\mathbb{Z}[G/G]$.

The only irreducibles occuring in $\mathbb{C}[G/B]$ are $U_1$ and $V_1$. Now, $\lambda_\chi(\epsilon_B) = 0$ only for $\chi = V_1$ by Lemma 6.14. On the other hand, $\lambda_\chi(\delta_B) = \lambda_\chi(|G|^2 (1 - \text{pr}_G) BN \times NB) = |G|^2 p^2 + p - 1 \neq 0$ for $\chi = V_1$ so again by the corollary, sequence 11 is exact at $\mathbb{Z}[G/B]$. □

Calculating the eigenvalues for $N'N \times NN'$ or $NN' \times N'N$ is more subtle. An indication of this can be seen by attempting to calculate the character sum

$$\sum_{n \in N} \sum_{n' \in N'} \overline{\chi(nn')} = \sum_{[g] \in C(G)} \overline{\chi([g])} c([g])$$

by summing over conjugacy classes as we did in Lemma 6.13. Here, $c([g])$ denotes the number of elements of the form $nn'$ which lie in $[g]$. It can be shown that $c([g])$ is essentially the number of points on an elliptic curve $E_{[g]}/\mathbb{F}_p$ which varies in a non-trivial way with $[g]$. Moreover, the Hasse bound for the number of points in $E_{[g]}(\mathbb{F}_p)$ is not sufficient to verify the hypotheses of Corollary 6.8.

## 7. A double coset algebra

In this section, we will describe the double coset algebra for $\mathbb{Z}[N\backslash G/N]$ explicitly as a $\mathbb{Z}$-module. This will be used later to get a handle on the operator $NN' \times N'N \in \mathbb{Z}[N\backslash G/N]$.

**Proposition 7.1.**

$$\mathbb{Z}[N\backslash G/N] = \mathbb{Z}N \oplus \bigoplus_{t \in \mathbb{F}_p^1/\sim} \mathbb{Z}N\sigma_t N$$

*where*

$$\sigma_t = \begin{pmatrix} 1 & 1 \\ 1 & t \end{pmatrix},$$

$\mathbb{F}_p^1 = \mathbb{F}_p - \{1\}$, *and* $t \sim t^{-1}$ *if* $t \neq 0$.

*Proof.* The entries of a matrix in $G$ will be denoted by $a, b, c, d$ starting in the upper right hand corner going clockwise.

Suppose $g = \in G - N$ is given so either $ac \neq 0$ or $bd \neq 0$. It will be shown that there exist $n_1, n_2 \in N$ such that $n_1 g n_2$ is one of $\sigma_t$ where $t \in \mathbb{F}_p^1$

At the expense of multiplication on the right or left by $\omega$, we obtain a matrix such that $ac \neq 0$ and $b \neq 0$. Then multiplication on the right and left respectively by

$$\begin{pmatrix} 1 & 0 \\ 0 & ac^{-1} \end{pmatrix}, \quad \begin{pmatrix} a^{-1} & 0 \\ 0 & b^{-1} \end{pmatrix},$$

for instance, yields the matrix $\sigma_{ad(bc)^{-1}}$. □

**Lemma 7.2.**
$$\deg(N) = 1$$
$$\deg(N\sigma_0 N) = 2(p-1)$$
$$\deg(N\sigma_{-1}N) = (p-1)/2$$
$$\deg(N\sigma_t N) = (p-1) \text{ for } t \in \mathbb{F}_p - \{-1, 0, 1\}.$$

*Proof.* Clearly, $\deg(N) = 1$. For $g = \sigma_t$ where $t \neq 1 \in \mathbb{F}_p$, it is easy to check that

$$\sigma_t C \sigma_t^{-1} \cap C = \left\{ \begin{pmatrix} x & 0 \\ 0 & x \end{pmatrix} \mid x \in \mathbb{F}_p^\times \right\}$$

$$\sigma_t C \sigma_t^{-1} \cap \omega C = \begin{cases} \left\{ \begin{pmatrix} x & 0 \\ 0 & -x \end{pmatrix} \mid x \in \mathbb{F}_p^\times \right\} & \text{if } t = -1 \\ \emptyset & \text{if } t \neq -1 \end{cases}$$

$$\sigma_t \omega C \sigma_t^{-1} \cap C = \begin{cases} \left\{ \begin{pmatrix} yt & 0 \\ 0 & y \end{pmatrix} \mid y \in \mathbb{F}_p^\times \right\} & \text{if } t \neq 0 \\ \emptyset & \text{if } t = 0 \end{cases}$$

$$\sigma_t \omega C \sigma_t^{-1} \cap \omega C = \begin{cases} \left\{ \begin{pmatrix} x & 0 \\ 0 & x \end{pmatrix} \mid x \in \mathbb{F}_p^\times \right\} & \text{if } t = -1 \\ \emptyset & \text{if } t \neq -1 \end{cases}$$

Hence, we have
$$|N_{\sigma_0}| = (p-1)$$
$$|N_{\sigma_{-1}}| = 4(p-1)$$
$$|N_{\sigma_t}| = 2(p-1) \text{ for } t \in \mathbb{F}_p - \{-1, 0, 1\}$$

The lemma follows by the definition of degree and the fact that $\#N = 2(p-1)^2$. □

## 8. A RELATION BETWEEN OPERATORS

In this section, we describe a relation between operators in $\mathbb{Z}[N\backslash G/N]$ which allows us to obtain information about the operator $NN' \times N'N$.

To describe this relation, it is necessary to introduce some more subgroups of $G$. Let $C''$ and $N''$ be the stabilisers in $G$ of $(1, -1)$ and $\{1, -1\}$, respectively. Then

$N'' = C'' \cup \omega''C''$ is the normaliser in $G$ of $C''$, where

$$\omega'' = \begin{pmatrix} 1 & 0 \\ 0 & -1 \end{pmatrix}.$$

The degrees of some operators associated to $N''$ are given for later reference.

**Lemma 8.1.**

$$\deg(NN'') = (p-1)/2$$
$$\deg(N''N) = (p-1)/2$$

We now decompose in terms of the natural basis given in Lemma 7.1 various operators in $\mathbb{Z}[N\backslash G/N]$ which are associated to the subgroups $N'$, $N''$, $B$, $G$.

**Lemma 8.2.**

$$NN' \times N'N = \frac{p-1}{2}N + \sum_{t \in \mathbb{F}_p^1/\sim, \left(\frac{t}{p}\right)=-1} N\sigma_t N.$$

*Proof.* We have that

$$NN' \times N'N = \frac{\deg(N'N)}{|N'|} \sum_{n' \in N'} \frac{\deg(NN')}{\deg(Nn'N)} Nn'N$$
$$= \frac{1}{8} \sum_{n' \in N'} \frac{1}{\deg(Nn'N)} Nn'N$$
$$= \frac{1}{4} \sum_{c' \in C'} \frac{1}{\deg(Nc'N)} Nc'N$$

where that last step follows from the fact that the involution $\omega'$ of $N'$ lies in $N$.

Consider the element

$$c' = \begin{pmatrix} x & \lambda y \\ y & x \end{pmatrix} \in C'.$$

If $xy = 0$, then $Nc'N = N$. Otherwise, $Nc'N = N\,\sigma_t N$ where $t = \frac{x^2}{\lambda y^2} \in \mathbb{F}_p^1$ is a non-square.

It is easy to check that there are $2(p-1)$ elements $c' \in C'$ such that $xy = 0$. Given $t \in \mathbb{F}_p^1$ a non-square, there are $2(p-1)$ elements $c' \in C'$ such that $t = x^2/\lambda y^2$. Since $N\sigma_t N = N\sigma_{t^{-1}} N$, for $t \in \mathbb{F}_p^1$ a non-square, there are $4(p-1)$ elements $c' \in C'$ such that $Nc'N = N\sigma_t N$ if $t \neq -1$ and $2(p-1)$ if $t = -1$. The result follows from the calculation of degrees in Lemma 7.2. $\square$

**Lemma 8.3.**

$$NN'' \times N''N = \frac{p-1}{2}N + \sum_{t \in \mathbb{F}_p^1/\sim, \left(\frac{t}{p}\right)=1} N\sigma_t N.$$

*Proof.* We have that

$$NN'' \times N''N = \frac{\deg(N''N)}{|N''|} \sum_{n'' \in N''} \frac{\deg(NN'')}{\deg(Nn''N)} Nn''N$$

$$= \frac{1}{8} \sum_{n'' \in N''} \frac{1}{\deg(Nn''N)} Nn''N$$

$$= \frac{1}{4} \sum_{c'' \in C''} \frac{1}{\deg(Nc''N)} Nc''N$$

where that last step follows from the fact that the involution $\omega''$ of $N''$ lies in $N$.

Consider the element

$$c'' = \begin{pmatrix} x & y \\ y & x \end{pmatrix} \in C''.$$

If $xy = 0$, then $Nc''N = N$. Otherwise, $Nc''N = N\,\sigma_t N$ where $t = \frac{x^2}{y^2} \in \mathbb{F}_p^1$ is a non-zero square.

It is easy to check that there are $2(p-1)$ elements $c'' \in C''$ such that $xy = 0$. Given $t \in \mathbb{F}_p^1$ a non-zero square, there are $2(p-1)$ elements $c'' \in C''$ such that $t = x^2/y^2$. Since $N\sigma_t N = N\sigma_{t^{-1}} N$, for $t \in \mathbb{F}_p^1$ a non-zero square, there are $4(p-1)$ elements $c'' \in C''$ such that $Nc''N = N\sigma_t N$ if $t \ne -1$ and $2(p-1)$ if $t = -1$. The result follows from the calculation of degrees in Lemma 7.2. $\square$

**Lemma 8.4.**

$$NB \times BN = 2N + N\sigma_0 N$$

*Proof.* We have that

$$NB \times BN = \frac{\deg(BN)}{|B|} \sum_{b \in B} \frac{\deg(NB)}{\deg(NbN)} NbN$$

$$= \frac{1}{|C|} \sum_{b \in B} \frac{2}{\deg(NbN)} NbN.$$

Consider the element

$$b = \begin{pmatrix} x & z \\ 0 & y \end{pmatrix}.$$

If $z = 0$, then $NbN = N$. Otherwise, we have $NbN = N\sigma_0 N$. There are $(p-1)^2$ elements $b \in B$ such that $z = 0$ and there are $(p-1)^3$ elements $b \in B$ such that $z \ne 0$. The result follows from the calculation of degrees in Lemma 7.2. $\square$

**Lemma 8.5.**

$$NG \times GN = N + \sum_{t \in \mathbb{F}_p^1/\sim} N\sigma_t N.$$

*Proof.* We have that

$$NG \times GN = \frac{\deg(GN)}{|G|} \sum_{g \in G} \frac{\deg(NG)}{\deg(NgN)} NgN$$

$$= \frac{1}{|N|} \sum_{g \in G} \frac{1}{\deg(NgN)} NgN$$

$$= \frac{1}{|N|}\left(|N|\, N + \sum_{t \in \mathbb{F}_p^1/\sim} \frac{1}{\deg(N\sigma_t N)} |N\sigma_t N|\, N\sigma_t N\right)$$

$$= N + \sum_{t \in \mathbb{F}_p^1/\sim} N\sigma_t N$$

where the last step follows from the fact that $|N\sigma_t N| = \deg(N\sigma_t N) \cdot |N|$. The result follows from Lemma 7.2. $\square$

The above lemmas yield the following relations.

**Proposition 8.6.** *We have the following relation of double coset operators*

$$(NN' \times N'N + NN'' \times N''N) + NB \times BN$$
$$= pN + NG \times GN.$$

*Proof.* The result follows from the preceeding lemmas. $\square$

**Proposition 8.7.**

$$N\sigma_{-1}N \times N\sigma_{-1}N = NN'' \times N''N$$

*Proof.* We have that

$$N\sigma_{-1}N \times N\sigma_{-1}N = \frac{\deg(N\sigma_{-1}N)}{|N|} \sum_{n \in N} \frac{\deg(N\sigma_{-1}N)}{\deg(N\sigma_{-1}n\sigma_{-1}N)} N\sigma_{-1}n\sigma_{-1}N$$

$$= \frac{\deg(N\sigma_{-1}N)}{|N|} \sum_{n'' \in N''} \frac{\deg(N\sigma_{-1}N)}{\deg(Nn''N)} Nn''N$$

$$= NN'' \times N''N$$

where the second last step follows from the fact that $\sigma_{-1}N\sigma_{-1} = N''$, and the last step follows from the fact that

$$\deg(N\sigma_{-1}N) = \deg(N''N) \text{ and } |N| = |N''|.$$

$\square$

## 9. Eigenvalues of a double coset operator

By the previous section, to know how $NN' \times N'N$ acts on $\mathbb{Z}[G/N]$, it suffices to know how $N\sigma_{-1}N$ acts on $\mathbb{Z}[G/N]$ since we know how the other operators in the relation act. In this section, we derive a formula for the eigenvalues of $N\sigma_{-1}N$.

**Lemma 9.1.** *Let $\chi$ be an irreducible character of $G$. Then the trace of $N\sigma N$ on the $\chi$-component of $\mathbb{C}[G/N]$ is given by*

$$\operatorname{tr}_\chi(N\sigma N) = \chi(1)\frac{\deg(N\sigma N)}{|N|} \sum_{g \in \sigma N} \overline{\chi(g)}.$$

*Proof.* Note that
$$\mathrm{tr}_\chi(N\sigma N) = \mathrm{tr}(\mathrm{pr}_\chi \circ N\sigma N).$$

Choose a set of inequivalent representatives $g_1, \ldots, g_n$ for $G/N$, where $n = |G/N|$. Then $g_1 N, \ldots, g_n N$ forms a $\mathbb{C}$-basis for $\mathbb{C}[G/N]$. To calculate the trace of the map $\mathrm{pr}_\chi \circ N\sigma N$, it suffices to compute for each $i$, the coefficient $\alpha_i$ of $g_i N$ in $(\mathrm{pr}_\chi \circ N\sigma N)(g_i N)$. The trace is then given by $\sum_{i=1}^n \alpha_i$.

To begin with, we have
$$(\mathrm{pr}_\chi \circ N\sigma N)(g_i N) = \left(\frac{\chi(1)}{|G|} \sum_{g \in G} \overline{\chi(g)} g\right)\left(\frac{\deg(N\sigma N)}{|N|} \sum_{n \in N} g_i n\sigma N\right)$$
$$= \frac{\chi(1)}{|G|} \frac{\deg(N\sigma N)}{|N|} \sum_{n \in N} \sum_{g \in G} \overline{\chi(g)} g_i n\sigma g N$$

where the last equality follows from the fact that $\mathrm{pr}_\chi$ is a $\mathbb{C}[G]$-module homomorphism. The coefficient of $g_i N$ in the above element of $\mathbb{C}[G/N]$ is then given by
$$\alpha_i = \frac{\chi(1)}{|G|} \frac{\deg(N\sigma N)}{|N|} \sum_{n \in N} \sum_{g \in G, g_i n\sigma g N = g_i N} \overline{\chi(g)}$$
$$= \frac{\chi(1)}{|G|} \frac{\deg(N\sigma N)}{|N|} \sum_{n \in N} \sum_{g \in \sigma^{-1} n^{-1} N} \overline{\chi(g)}$$
$$= \frac{\chi(1)}{|G|} \frac{\deg(N\sigma N)}{|N|} |N| \sum_{g \in \sigma N} \overline{\chi(g)}$$

since $\sigma^{-1} N = \sigma N$ for the representatives $\sigma$ we take for $N\sigma N$. □

**Lemma 9.2.** *Let $g \in G$ be a non-scalar element with trace $t$ and determinant $n$. The number of elements in $\sigma_{-1} N$ which are conjugate to $g$ is equal to*
$$c([g]) = 2\left(1 + \left(\frac{t^2 - 2n}{p}\right)\right).$$
*where $[g]$ denotes the conjugacy class of $g$.*

*Proof.* The general element in $\sigma_{-1} C$
$$\begin{pmatrix} 1 & 1 \\ 1 & -1 \end{pmatrix} \begin{pmatrix} x & 0 \\ 0 & y \end{pmatrix}$$
has trace $t = x - y$ and determinant $n = -2xy$. Therefore, $x$ satisfies $2x^2 - 2tx + n = 0$. Conversely, if $x$ satisfies $2x^2 - 2tx + n = 0$, then putting $y = x - t$ yields an element $h \in \sigma_{-1} C$ with trace $t$ and determinant $n$.

The general element in $\sigma_{-1} \omega C$
$$\begin{pmatrix} 1 & 1 \\ 1 & -1 \end{pmatrix} \begin{pmatrix} 0 & x \\ y & 0 \end{pmatrix}$$
has trace $t = x + y$ and determinant $n = 2xy$. Therefore, $x$ also satisfies $2x^2 - 2tx + n = 0$, and conversely if $x$ satisfies $2x^2 - 2tx + n = 0$, then putting $y = t - x$ yields an element $h \in \sigma_{-1}\omega C$ with trace $t$ and determinant $n$.

The number of solutions to $2x^2 - 2tx + n = 0$ is $\left(1 + \left(\frac{t^2 - 2n}{p}\right)\right)$ so $c(g) = 2(1 + \left(\frac{t^2-2n}{p}\right))$. □

**Proposition 9.3.** *Let $\chi$ be the character of $W_{\alpha,\beta}$ where $\beta = \alpha^{-1}$, $\alpha^{\frac{p-1}{2}} = 1$ and $\alpha^2 \neq 1$ so that $\alpha(-1) = 1$. Then*

$$\operatorname{tr}_\chi(N\sigma_{-1}N) = \frac{p+1}{2} \sum_{d \in \mathbb{F}_p^\times} \alpha(d)\left(\frac{1+d^2}{p}\right)$$

*and hence*

$$\lambda_\chi(N\sigma_{-1}N) = \frac{1}{2} \sum_{d \in \mathbb{F}_p^\times} \alpha(d)\left(\frac{1+d^2}{p}\right).$$

*Proof.* Note that $\overline{\chi} = \chi$. We compute

$$\sum_{g \in \sigma_{-1}N} \chi(g) = \sum_{[g] \in C(G)} \chi([g]) \cdot c([g])$$

where $C(G)$ is the set of conjugacy classes in $G$. Recall the description of $C(G)$ in Section 5. We first note that there are no scalars $a_x$ in $\sigma_{-1}N$. There are precisely $2(1+\left(\frac{2}{p}\right))$ elements in $\sigma_{-1}N$ which are conjugate to a given $[b_x]$. Finally, there are $2(1+\left(\frac{x^2+y^2}{p}\right))$ elements in $\sigma_{-1}N$ which are conjugate to a given $[\kappa_{x,y}]$.

From the values of $\chi$ on each type of conjugacy class (see Table 1), we obtain that

$$\sum_{[g] \in C(G)} \chi([g]) \cdot c([g]) = 2(p-1) \sum_{d \in \mathbb{F}_p^\times} = \alpha(d)\left(\frac{1+d^2}{p}\right)$$

where $d = y/x$ and $d \sim d^{-1}$. The result follows from the fact that $\chi(1)\frac{\deg(N\sigma_{-1}N)}{|N|} = \frac{p+1}{4(p-1)}$. □

**Proposition 9.4.** *Let $\chi$ be the character of $X_\phi$ where $\phi^{p+1} = 1$, $\phi^{\frac{p+1}{2}} \neq 1$, and $\phi^{p-1} \neq 1$ so that $\phi\mid_{\mathbb{F}_p^\times} = 1$. Then*

$$\operatorname{tr}_\chi(N\sigma_{-1}N) = -\frac{1}{2} \sum_{\gamma \in C'} \phi(\gamma)\left(\frac{\gamma^2 + \overline{\gamma}^2}{p}\right)$$

*and hence*

$$\lambda_\chi(N\sigma_{-1}N) = -\frac{1}{2(p-1)} \sum_{\gamma \in C'} \phi(\gamma)\left(\frac{\gamma^2 + \overline{\gamma}^2}{p}\right).$$

*Proof.* Note that $\overline{\chi} = \chi$. We compute

$$\sum_{g \in \sigma_{-1}N} \chi(g) = \sum_{[g] \in C(G)} \chi([g]) \cdot c([g])$$

where $C(G)$ is the set of conjugacy classes in $G$.

We first note that there are no scalars $a_x$ in $\sigma_{-1}N$. There are precisely $2(1+\left(\frac{2}{p}\right))$ elements in $\sigma_{-1}N$ which are conjugate to a given $[b_x]$. Finally, there are $2(1+\left(\frac{\gamma^2+\overline{\gamma}^2}{p}\right))$ elements in $\sigma_{-1}N$ which are conjugate to a given $[\gamma_{x,y}]$.

From the values of $\chi$ on each type of conjugacy class, we obtain that

$$\sum_{[g] \in C(G)} \chi([g]) \cdot c([g]) = -2 \sum_{\gamma \in C'} \phi(\gamma)\left(\frac{\gamma^2 + \overline{\gamma}^2}{p}\right)$$

The result follows from the fact that $\chi(1)\frac{\deg(N\sigma_{-1}N)}{|N|} = \frac{1}{4}$. □

**Proposition 9.5.** *Suppose that $p \equiv 1 \pmod{4}$. Let $\chi$ be the character of $V_\alpha$ where $\alpha = \left(\frac{\cdot}{p}\right)$. Then*

$$\operatorname{tr}_\chi(N\sigma_{-1}N) = \frac{p}{4}\sum_{d\in\mathbb{F}_p^\times} \alpha(d)\left(\frac{1+d^2}{p}\right) - \frac{p}{4(p-1)}\sum_{\gamma\in C'} \alpha(\gamma^{p+1})\left(\frac{\gamma^2+\overline{\gamma}^2}{p}\right)$$

*and hence*

$$\lambda_\chi(N\sigma_{-1}N) = \frac{1}{4}\sum_{d\in\mathbb{F}_p^\times} \alpha(d)\left(\frac{1+d^2}{p}\right) - \frac{1}{4(p-1)}\sum_{\gamma\in C'} \alpha(\gamma^{p+1})\left(\frac{\gamma^2+\overline{\gamma}^2}{p}\right).$$

*Proof.* Note that $\overline{\chi} = \chi$. We compute

$$\sum_{g\in\sigma_{-1}N} \chi(g) = \sum_{[g]\in C(G)} \chi([g]) \cdot c([g])$$

where $C(G)$ is the set of conjugacy classes in $G$. Recall the description of $C(G)$ in Section 5. We first note that there are no scalars $a_x$ in $\sigma_{-1}N$. There are $2(1+\left(\frac{x^2+y^2}{p}\right))$ elements in $\sigma_{-1}N$ which are conjugate to a given $[\kappa_{x,y}]$. Finally, there are $2(1+\left(\frac{\gamma^2+\overline{\gamma}^2}{p}\right))$ elements in $\sigma_{-1}N$ which are conjugate to a given $[\gamma_{x,y}]$.

From the values of $\chi$ on each type of conjugacy class (see Table 1), we obtain that

$$\sum_{[g]\in C(G)} \chi([g]) \cdot c([g])$$
$$= \sum_{(x,y)\in\mathbb{F}_p^\times\times\mathbb{F}_p^\times-\Delta/\sim} \alpha(xy) \cdot 2(1+\left(\frac{x^2+y^2}{p}\right))$$
$$- \sum_{\gamma\in C'-\mathbb{F}_p^\times/\sim} \alpha(\gamma^{p+1}) \cdot 2(1+\left(\frac{\gamma^2+\overline{\gamma}^2}{p}\right))$$

We compute each sum separately.

$$\sum_{(x,y)\in\mathbb{F}_p^\times\times\mathbb{F}_p^\times-\Delta/\sim} \alpha(xy) \cdot 2(1+\left(\frac{x^2+y^2}{p}\right)) = -(p-1)\cdot(1+\left(\frac{2}{p}\right))$$
$$+ (p-1)\sum_{d\in\mathbb{F}_p^\times} \alpha(d)\cdot\left(\frac{1+d^2}{p}\right)$$

where $d = y/x$ and $d \sim d^{-1}$.

$$\sum_{\gamma\in C'-\mathbb{F}_p^\times/\sim} \alpha(\gamma^{p+1}) \cdot 2(1+\left(\frac{\gamma^2+\overline{\gamma}^2}{p}\right)) = -(p-1)\cdot(1+\left(\frac{2}{p}\right))$$
$$+ \sum_{\gamma\in C'} \alpha(\gamma^{p+1})\cdot\left(\frac{\gamma^2+\overline{\gamma}^2}{p}\right)$$

The result follows from the fact that $\chi(1)\frac{\deg(N\sigma_{-1}N)}{|N|} = \frac{p}{4(p-1)}$. □

**Remark 9.6.** *The character sum in Proposition 9.3 is an instance of a Legendre character sum [27]. The character sum in Proposition 9.4 is an instance of a Soto-Andrade sum [18] [30].*

## 10. Non-vanishing of $\lambda_\chi(NN' \times N'N)$

We show the non-vanishing of the eigenvalues $\lambda_\chi(NN' \times N'N)$ for $\chi$ occuring in $\mathbb{C}[G/N']$ which will be used in the next section to conclude that sequence 11 is almost-exact.

Let $\chi$ be the character of an irreducible component of $\mathbb{C}[G/N']$ which is not trivial. By the relation in Proposition 8.6, we see that

$$NN' \times N'N = pN - NN'' \times N''N - NB \times BN + NG \times GN.$$

Now, the eigenvalues of $NB \times BN$ and $NG \times GN$ acting on the $\chi$-component of $\mathbb{C}[G/N]$ must be zero because these $\mathbb{C}[G]$-module homomorphisms factor through $\mathbb{C}[G/B]$ and $\mathbb{C}[G/G]$ which do not contain the irreducible representation $\chi$. Hence,

$$\lambda_\chi(NN' \times N'N) = p - \lambda_\chi(NN'' \times N''N).$$

On the other hand, by the relation in Proposition 8.7, we see that

$$\lambda_\chi(NN'' \times N''N) = \lambda_\chi(N\sigma_{-1}N)^2.$$

Thus, to show that $\lambda_\chi(NN' \times N'N) \neq 0$ is equivalent to showing that $\lambda_\chi(N\sigma_{-1}N) \neq \pm\sqrt{p}$.

**Proposition 10.1.** *Let $\chi$ be the character of $W_{\alpha,\beta}$ where $\beta = \alpha^{-1}$, $\alpha^{\frac{p-1}{2}} = 1$ and $\alpha^2 \neq 1$ so that $\alpha(-1) = 1$. Then*

$$\lambda_\chi(N\sigma_{-1}N) \neq \pm\sqrt{p}$$

*Proof.* From Proposition 9.3, we see that $\lambda_\chi(N\sigma_{-1}N) \in \mathbb{Q}(\zeta_{p-1})$ where $\zeta_{p-1}$ is a primitive $p-1$-th root of unity. But this field does not contain $\mathbb{Q}(\sqrt{p})$ as $p$ does not ramify in it. □

**Proposition 10.2.** *Let $\chi$ be the character of $X_\phi$ where $\phi^{p+1} = 1$, $\phi^{\frac{p+1}{2}} \neq 1$, and $\phi^{p-1} \neq 1$ so that $\phi|_{\mathbb{F}_p^\times} = 1$. Then*

$$\lambda_\chi(N\sigma_{-1}N) \neq \pm\sqrt{p}$$

*Proof.* From Proposition 9.4, we see that $\lambda_\chi(N\sigma_{-1}N) \in \mathbb{Q}(\zeta_{p+1})$ where $\zeta_{p+1}$ is a primitive $p+1$-th root of unity. But this field does not contain $\mathbb{Q}(\sqrt{p})$ as $p$ does not ramify in it. □

**Proposition 10.3.** *Suppose $p \equiv 1 \pmod{4}$. Let $\chi$ be the character of $V_\alpha$ where $\alpha = \left(\frac{\cdot}{p}\right)$. Then*

$$\lambda_\chi(N\sigma_{-1}N) \neq \pm\sqrt{p}$$

*Proof.* From Proposition 9.5, we see that $\lambda_\chi(N\sigma_{-1}N) \in \mathbb{Q}$ which does not contain $\pm\sqrt{p}$. □

## 11. Exactness at $\mathbb{C}[G/N]$ and $\mathbb{C}[G/N']$

In this section, we show almost-exactness at $\mathbb{Z}[G/N]$ and $\mathbb{Z}[G/N']$ in sequence 11 using Corollary 6.8 and non-vanishing results of the previous section.

**Proposition 11.1.** *Let $\chi$ be the character of the trivial representation $U_1$. Then*
$$\lambda_\chi(NN' \times N'N) = \text{tr}_\chi(NN' \times N'N) = \frac{p^2 - 1}{4}.$$

*Proof.* We could use Lemma 6.11 directly, but let's amuse ourselves and compute it using Lemmas 8.2 and 9.1.

$$\begin{aligned}
\lambda_\chi &= \text{tr}_\chi(NN' \times N'N) \\
&= \frac{p-1}{2} \text{tr}_\chi(N) + \sum_{t \in \mathbb{F}_p^1/\sim, \left(\frac{t}{p}\right)=-1} \text{tr}_\chi(N\sigma_t N) \\
&= \frac{\deg(N)}{|N|} \cdot |N| + \sum_{t \in \mathbb{F}_p^1/\sim, \left(\frac{t}{p}\right)=-1} \frac{\deg(N\sigma_t N)}{|N|} \cdot |N| \\
&= \frac{p^2 - 1}{4}
\end{aligned}$$

where the last step follows from the fact that the sum is $\frac{p-1}{4} \cdot (p-1) = \frac{(p-1)^2}{4}$ if $\left(\frac{-1}{p}\right) = -1$ and $\frac{\frac{p-1}{2}-1}{2} \cdot (p-1) + \frac{p-1}{2} = \frac{(p-1)^2}{4}$ if $\left(\frac{-1}{p}\right) = 1$. $\square$

**Theorem 5.** *Let $\chi$ be the character of $W_{\alpha,\alpha^{-1}}$, $X_\phi$, or $V_{\left(\frac{\cdot}{p}\right)}$ if $p \equiv 1 \pmod{4}$ as in Propositions 9.3, 9.4, and 9.5. Then $\lambda_\chi(NN' \times N'N) \neq 0$.*

*Proof.* By Propositions 10.1, 10.2, 10.3, $\lambda_\chi(N\sigma_{-1}N) \neq \pm\sqrt{p}$. The result then follows from the relations in Propositions 8.6 and 8.7. $\square$

**Corollary 11.2.** *Sequence 11 is almost-exact at $\mathbb{Z}[G/N]$ and $\mathbb{Z}[G/N']$.*

*Proof.* By Theorem 5 and Proposition 11.1, $\lambda_\chi(\epsilon_N) = \lambda_\chi(NN' \times N'N) = \lambda_\chi(N'N \times N'N)$ is non-zero for each irreducible component of $\mathbb{C}[G/N']$. Thus, by Corollary 6.8, sequence 11 is almost-exact at $\mathbb{Z}[G/N']$.

By Lemma 6.13, $\lambda_\chi(\epsilon_N) = \lambda_\chi(|G|^2 (1 - \text{pr}_G)NB \times BN) = 0$ for all irreducible components of $\mathbb{C}[G/N]$ except $V_1$. However, by Lemma 5.4, these are precisely the irreducible components occuring in $\mathbb{C}[G/N']$. Furthermore, by Theorem 5 and Proposition 11.1, $\lambda_\chi(\delta_N) = \lambda_\chi(NN' \times N'N) \neq 0$ for such $\chi$. Applying Corollary 6.8 again shows that sequence 11 is almost-exact at $\mathbb{Z}[G/N]$. $\square$

Combining Corollary 11.2 and Proposition 6.15, and then applying Proposition 3.7, we obtain Theorem 2.

## 12. Some numerical data

We have computed for primes $p \leq 19$, the determinant of the map $N'N \times NN'$ on $\mathbb{C}[G/N']$ and on its components $U = U_1$, $W = \oplus_\alpha W_{\alpha,\alpha^{-1}}$, $X = \oplus_\phi X_\phi$, and $V = V_{\left(\frac{\cdot}{p}\right)}$ if $p \equiv 1 \pmod{4}$, using the expressions in terms of character sums obtained in previous sections.

For the primes $p \leq 13$, the computation of the determinant on $\mathbb{C}[G/N']$ was done independently by computing the matrix of $N'N \times NN'$ and then calculating its determinant. This value coincides with the one obtained by evaluating character sums.

TABLE 2. Determinants of $N'N \times NN'$

| $p$ | $|G|$ | $\mathbb{C}[G/N']$ | $U$ | $W$ | $X$ | $V$ |
|---|---|---|---|---|---|---|
| 3 | $2^4 \cdot 3$ | $2^3$ | 2 | 1 | 4 | |
| 5 | $2^5 \cdot 3 \cdot 5$ | $2^{11} \cdot 3$ | $2 \cdot 3$ | 1 | 1 | $2^{10}$ |
| 7 | $2^5 \cdot 3^2 \cdot 7$ | $2^{20} \cdot 3^9$ | $2^2 \cdot 3$ | $3^8$ | $2^{18}$ | |
| 11 | $2^4 \cdot 3 \cdot 5^2 \cdot 11$ | $2^{71} \cdot 3 \cdot 5^{13}$ | $2 \cdot 3 \cdot 5$ | $5^{12}$ | $2^{70}$ | |
| 13 | $2^5 \cdot 3^2 \cdot 7 \cdot 13$ | $2^{83} \cdot 3^{15} \cdot 7$ | $2 \cdot 3 \cdot 7$ | $2^{56} \cdot 3^{14}$ | 1 | $2^{26}$ |
| 17 | $2^9 \cdot 3^2 \cdot 17$ | $2^{215} \cdot 3^2 \cdot 19^{32}$ | $2^3 \cdot 3^2$ | $2^{144}$ | $19^{32}$ | $2^{68}$ |
| 19 | $2^4 \cdot 3^4 \cdot 5 \cdot 19$ | $2^{163} \cdot 3^{78} \cdot 5 \cdot 17^{40}$ | $2 \cdot 3^2 \cdot 5$ | $3^{40} \cdot 17^{40}$ | $2^{162} \cdot 3^{36}$ | |

## 13. CONCLUSIONS

There are several questions still remaining. First of all, it would be interesting to calculate in detail the homology of sequence 11 and ultimately the one in Theorem 1. A related but perhaps different question concerns describing a minimal value for the degree of the isogenies which may occur in Theorem 1 in some suitable sense.

It would be interesting to explain the relations of double coset operators in Propositions 8.6 and 8.7, in particular why the normaliser of a split Cartan subgroup $N''$ which is conjugate to the standard one arises.

One might consider how Merel's conjecture generalises to other relations of jacobians. For instance, there is a variant of the relation in Theorem 1 which was also proved in [4] [5], and [10].

**Theorem 6** (Chen,Edixhoven).

$$J_{C'} \times J_B{}^2 \text{ is isogenous over } \mathbb{Q} \text{ to } J_C \times J_G{}^2$$

In this case already, it is not clear which correspondences give the desired relation of jacobians, as the canonical choices do not appear to work (thanks are due to L. Merel for pointing this out). A perhaps related problem is that the associated induced representations do not have multiplicity one.

## REFERENCES

[1] R.D. Accola. Riemann surfaces with automorphism groups admitting partitions. *Proceedings of the American Mathematical Society*, 21:477–482, 1969.

[2] R.D. Accola. Two theorems on riemann surfaces with non-cyclic automorphism groups. *Proceedings of the American Mathematical Society*, 25:598–602, 1970.

[3] N. Arsenault. The class groups of arithmetically equivalent algebras. Master's thesis, McGill University, 1996.

[4] I. Chen. *The Jacobian of Modular Curves Associated to Cartan Subgroups*. PhD thesis, University of Oxford, 1996.

[5] I. Chen. The Jacobians of non-split Cartan modular curves. *Proceedings of the London Mathematical Society*, 77, Part 1:1–38, 1998.

[6] C. Curtis and I. Reiner. *Methods of Representation Theory*. John Wiley & Sons, Inc., 1981.

[7] H. Darmon. The equations $x^n + y^n = z^2$ and $x^n + y^n = z^3$. *International Mathematics Research Notices*, 72(1):263–273, 1993.

[8] H. Darmon and L. Merel. Winding quotients and some variants of Fermat's last theorem. *Journal für die Reine und Angewandte Mathematik*, 490:81–100, 1997.

[9] P. Deligne and M. Rapoport. Les schémas de modules de courbes elliptiques. In P. Deligne and W. Kuyk, editors, *Modular Functions of One Variable II*, number 349 in Lecture Notes in Mathematics, pages 143–316. Springer-Verlag, 1972.

[10] S.J. Edixhoven. On a result of Imin Chen. Preprint (Duke algebraic geometry preprint server), 1996.

[11] R. Evans. Hermite character sums. *Pacific Journal of Mathematics*, 122(2):357–390, 1986.

[12] W. Fulton and J. Harris. *Representation Theory, A First Course*. Number 129 in Graduate Texts in Mathematics, Readings in Mathematics. Springer-Verlag, 1991.

[13] J. Greene. Hypergeometric functions over finite fields. *Transactions of the American Mathematical Society*, 301(1):77–101, 1987.

[14] J. Greene. Hypergeometric functions over finite fields and representations of SL(2, $q$). *Rocky Mountain Journal of Mathematics*, 23(2):547–568, 1993.

[15] B. Gross. *Arithmetic on Elliptic Curves with Complex Multiplication*. Number 776 in Lecture Notes in Mathematics. Springer-Verlag, 1980.

[16] A. Helversen-Pasotto. L'identité de Barnes for les corps finis. *Comptes Rendues de l'Académie des Sciences*, 286:297–300, 1978. Paris, Série A.

[17] E. Kani and M. Rosen. Idempotent relations and factors of jacobians. *Mathematische Annalen*, 284(2):307–327, 1989.

[18] N. Katz. Estimates for Soto-Andrade sums. *Journal für die Reine und Angewandte Mathematik*, 438:143–161, 1993.

[19] S. Lang. *Abelian varieties*. Number 7 in Interscience tracts in pure and applied mathematics. Interscience Publishers, Inc., New York, 1959.

[20] W. Li and J. Soto-Andrade. Barnes' identities and representations of GL(2), Part I: Finite field case. *Journal für die Reine und Angewandte Mathematik*, 344:171–197, 1983.

[21] G. Ligozat. Courbes modulaires de niveau 11. In J.P. Serre and D.B. Zagier, editors, *Modular Functions of One Variable V*, number 601 in Lecture Notes in Mathematics, pages 149–237. Springer-Verlag, 1977.

[22] B. Mazur. Modular curves and the Eisenstein ideal. *I.H.E.S. Publications Mathématiques*, 47:33–186, 1977.

[23] B. Mazur. Rational isogenies of prime degree. *Inventiones mathematicae*, 44:129–162, 1978.

[24] L. Merel. Arithmetic of elliptic curves and diophantine equations. preprint, November 18, 1996.

[25] F. Momose. Rational points on the modular curves $X_{\mathrm{split}}(p)$. *Compositio Mathematica*, 52:115–137, 1984.

[26] K. Ribet. On the equation $a^p + 2^\alpha b^p + c^p = 0$. *Acta Arithmetica*, 79(1):7–16, 1997.

[27] Y. Sawabe. Legendre character sums. *Hiroshima Mathematical Journal*, 22:15–22, 1992.

[28] J.P. Serre. Propriétés galoisiennes des points d'ordre fini des courbes elliptiques. *Inventiones Mathematicae*, 15:259–331, 1972.

[29] G. Shimura. *Introduction to the Arithmetic Theory of Automorphic Functions*. Iwanami Shoten, Publishers and Princeton University Press, 1971.

[30] J. Soto-Andrade. Geometrical Gelfand models, tensor quotients, and Weil representations. In *Proceedings of Symposia in Pure Mathematics*, volume 47, pages 305–316, 1987.

[31] J. Tate. On the conjectures of Birch and Swinnerton-dyer and a geometric analogue. Séminaire Bourbaki, no. 306, 1965-66.

[32] A. Weil. *Courbes algébriques et variétés abéliennes*. Hermann et Cie., Paris, 1971. Combined volume.

Department of Mathematics and Statistics, McGill University, Montreal, Quebec, Canada, H3A 2K6

*E-mail address*: `chen@math.mcgill.ca`